\documentclass[12pt,leqno]{article}
 
\def\w{\dot{w}}
\def\v{{\rm v}}
\def\sigmad{\dot{\sigma}} 
\def\int{\mathbb{Z}}

\def\Ue{{\cal U}_{\varepsilon}({\mathfrak g})}

\def\O{{\cal O}}

\def\proof{{\bf Proof. }}

\def\Pf{\proof}

\def\pf{\proof}

\def\V{{\cal V}}
\usepackage{times}
\usepackage{graphics}
\usepackage{amssymb}
\usepackage{amscd}
\usepackage{amsmath}
\input xy 
\xyoption{all}
\title{Spherical conjugacy classes and the Bruhat decomposition}
\newtheorem{theorem}{Theorem}[section]
\newtheorem{lemma}[theorem]{Lemma}

\newtheorem{proposition}[theorem]{Proposition}
\newtheorem{definition}[theorem]{Definition}
\newtheorem{remark}[theorem]{Remark}

\author{Giovanna Carnovale\\
Dipartimento di Matematica Pura ed Applicata\\
Torre Archimede - via Trieste 63 - 35121 Padova - Italy\\
email: carnoval@math.unipd.it }
\date{}
\begin{document}
\maketitle
\begin{abstract}Let $G$ be a connected, reductive algebraic group over an
  algebraically closed field of zero or good and odd characteristic. Let $B$ be a Borel subgroup of $G$. We characterize spherical conjugacy classes in $G$ as those intersecting only the double
  cosets $BwB$ in $G$ corresponding to involutions in the Weyl group of $G$.
\end{abstract}

\noindent{\bf Key words:} conjugacy class, spherical homogeneous space, Bruhat decomposition

\noindent{\bf MSC:} 20GXX; 20E45; 20F55; 14M15

\section*{Introduction}

The Bruhat decomposition of a connected reductive algebraic group $G$ over an algebraically closed field states that the two-sided cosets of $G$ with respect to a Borel subgroup $B$ (Bruhat cells) are naturally parametrized by the elements in the Weyl group of $G$ and have a well-understood geometrical behaviour. It is a
fundamental tool in the theory of algebraic groups, as it is 
relevant for the comprehension of the geometry of the flag
variety $G/B$, for instance, in the computation of its cohomology.  
Besides, intersection of Bruhat cells corresponding to opposite Borel subgroups (double Bruhat cells) play a significant role
in the description of the symplectic leaves of a natural Poisson
structure on $B$ (\cite{DCKP3}). 
New interest has been raised by Bruhat cells and
double Bruhat cells for their applications to total positivity (\cite{FZ}) and to
the theory of cluster algebras. For instance, as it has been very recently
shown, double Bruhat cells serve as a geometric model for cluster
algebras of finite type, since every cluster algebra of finite type
with principal coefficients at an arbitrary acyclic initial cluster can be
realized as the coordinate ring of a certain double Bruhat cell (\cite{YZ}).

The interplay between conjugacy classes in an algebraic group and the
Bruhat decomposition has been successfully exploited in the past. 
Probably the first results
in this sense are in \cite{regular} where the Bruhat decomposition of
a semisimple algebraic group $G$ is
used for the construction of a cross-section for the collection of
regular conjugacy classes of $G$. 

More recently, \cite{EG} and \cite{EG2} have provided an analysis of the intersection of conjugacy
classes in a Chevalley group with Bruhat cells corresponding to generalized Coxeter elements and
their conjugates.

If we consider spherical conjugacy classes, that is,
those conjugacy classes of a group $G$ on which $B$ acts
with finitely many orbits, it is natural to inquire about their intersection  with Bruhat cells.
A characterization of spherical conjugacy classes has been given in terms of a formula
involving the dimension of the class $\O$ and the maximal element $w$ in
the Weyl group $W$ of $G$ for which $\O\cap BwB$ is non-empty. This is
obtained in \cite{ccc} over the complex numbers and  in \cite{gio} over an arbitrary algebraically
closed field of zero or odd good characteristic. The motivation in \cite{ccc} was
the proof - in the spherical case - of a conjecture due to De Concini, Kac and Procesi on the
dimension of irreducible representations of quantum groups at the
roots of unity (\cite{DCKP}). The proof relied on the classifications of spherical
nilpotent orbits (\cite{pany2}) and of reductive spherical pairs (\cite{brion}) and on geometric properties
of spherical homogeneous spaces in the complex setting (\cite{brion},\cite{pany}). In \cite{gio}  a
different approach  was developed and a crucial step in the argument was
that every spherical conjugacy class intersects only Bruhat
cells $BwB$ for $w$ an involution in $W$. The aim of the present
paper is to show that this property fully characterizes spherical
conjugacy classes. 

\medskip

\noindent{\bf Theorem} Let $G$ be a connected reductive algebraic
group over an algebraically closed field of zero or good, odd characteristic.
A conjugacy class $\O$ in $G$ is spherical if and only if $\O$
intersects only Bruhat cells corresponding to involutions in the Weyl group of $G$.

\medskip

The paper is structured as follows: after fixing notation and
recalling basic facts about spherical homogeneous spaces and conjugacy classes in \S \ref{preliminaries},
we analyse the case of $G$ simple of type $G_2$ in full detail in \S \ref{sect-g2}. The reason for doing so is twofold. On the one hand we would like to give an idea of the techniques involved through an example, and on the other hand it would not be more efficient to treat the case of $G_2$ together with the others because separate descriptions for behaviour of roots with different length ratios are needed. 

In \S \ref{maximal} we restrict our attention to those conjugacy classes intersecting only Bruhat cells corresponding to involutions. For such a class $\O$ we consider the maximal element $w\in W$ for which $\O\cap BwB$ is non-empty and the
 set of $B$-orbits in $\O$ that are contained in $BwB$, the so-called maximal $B$-orbits. The properties of a special class of
representatives $x$ of maximal $B$-orbits are analyzed, allowing a description of
the centralizer $B_x$ in $B$. This is achieved by using the same
 strategy as in \cite{gio}. The proofs therein are rather laborious and
 need a case-by-case analysis but they apply also to the present situation so we use them referring to \cite{gio}. The hypothesis on the class $\O$ imposes restrictions on the
 representatives $x$ in maximal $B$-orbits: for instance, if
 $x=\w\v\in N(T)U$ then $\v$ lies in the subgroup generated by the root
 subgroups $X_\alpha$ for which $w\alpha=-\alpha$.  This condition is powerful for a general $w$ 
but it is empty when $w$ is the longest element $w_0$ in $W$ and   
%
it acts as $-1$ in the geometric representation. For this reason we deal with this situation separately and an unpleasant case-by-case analysis is needed in the doubly-laced case. This is done in \S \ref{bigcell}, where the theorem in this case is proved by showing the sufficient condition that the maximal $B$-orbits are finitely-many.  

The rest of the paper is devoted to an estimate of the centralizer $G_x$ in
$G$ of a representative $x$ in a maximal $B$-orbit. Indeed, since $\O$
is parted into finitely many $B$-orbits if and only if it has a dense
$B$-orbit (\cite{Bri,gross,knop,Vin}), we may conclude that $\O$ is
spherical once we prove that the dimension of a maximal $B$-orbit
equals the dimension of $\O$. In \S \ref{curves} we consider the
general case and we construct some families of elements contained in
$G_x\cap X_{\alpha}s_\alpha B$ for different roots $\alpha$. We need
different strategies according to the behaviour of $\alpha$ with respect to $w$. In particular,
when $w\alpha=-\alpha$ we apply the results in \S \ref{bigcell}.  
Once we have constructed enough elements in $G_x$ we show using the intersections $G_x\cap B\sigma B$ and induction on the length of $\sigma$ that the image
of $G_x$ through the projection of $G$  on $G/B$ is dense in the flag
variety obtaining the sought equality of dimensions.  
%
%

\section{Preliminaries}\label{preliminaries}

Unless otherwise stated $G$ will denote
a connected, reductive algebraic group over an algebraically closed field
$k$ of characteristic $0$ or odd and good (\cite[\S I.4]{131}). When we write an integer as an element in $k$ we shall mean the image of that integer in the prime field of $k$.

Let $B$ be a
Borel subgroup of $G$, let $T$ be a maximal torus contained in $B$ and let $B^-$ be the
Borel subgroup opposite to $B$. Let $U$
(respectively $U^-$) be the unipotent radical of $B$ (respectively
$B^-$). 

We shall denote by  $\Phi$ the root system 
relative to $(B, T)$; by
$\Delta=\{\alpha_1, \dots, \alpha_n\}$
the corresponding set of simple roots and by $\Phi^+$ the corresponding set of positive roots. We shall use the numbering of
the simple roots in \cite[Planches I-IX]{bourbaki}.

We shall denote by $W$ the
Weyl group associated with $G$ and by
$s_{\alpha}$ the reflection corresponding to the root $\alpha$. By
$\ell(w)$ we shall denote the length of the element $w\in W$ and
by ${\rm rk}(1-w)$ we shall mean the rank of $1-w$ in the
geometric representation of the Weyl group. By $w_0$ we shall denote
the  longest element in $W$ and $\vartheta$ will be the
automorphism of $\Phi$ given by $-w_0$. By $\Pi$ we shall always denote
a subset of $\Delta$ and $\Phi(\Pi)$ will indicate the corresponding
root subsystem of $\Phi$. We shall denote by
$W_\Pi$ the parabolic subgroup of $W$ generated by the $s_\alpha$ for $\alpha$ in $\Pi$. Given an element $w\in W$ we shall denote by $\w$
a representative of $w$ in the normalizer $N(T)$ of $T$.
For any root $\alpha$ in $\Phi$ we shall write $x_{\alpha}(t)$ for the
elements in the corresponding root subgroup $X_{\alpha}$ of
$G$. Moreover, we choose $x_\alpha(1)$ and $x_\alpha(-1)$ so that 
$x_\alpha(1)x_{-\alpha}(-1)x_{\alpha}(1)=n_\alpha\in s_\alpha T$ so
that the properties in \cite[Lemma 8.1.4]{springer} hold.

If $\Pi\subset\Delta$ we shall indicate by $P_\Pi$ the standard parabolic subgroup of $G$ whose Levi component contains the root subgroups corresponding to roots in $\Phi(\Pi)$ and by $P^u_\Pi$ its unipotent radical. If $\Pi=\{\alpha\}$ we shall simply write $P_\alpha$ and $P_\alpha^u$.

For $w\in W$, we will put
\begin{equation}\label{fivu}\Phi_w:=\{\alpha\in\Phi^+~|~w^{-1}\alpha\in-\Phi^+\}\end{equation}
\begin{equation}\label{uvu}U^w=\langle X_\alpha~|~\alpha\in\Phi_w\rangle,\quad U_w=\langle X_\alpha~|~\alpha\in\Phi^+\setminus\Phi_w\rangle
\end{equation}
so that $BwB=U^w\w B$ for any choice of $\w\in N(T)$.
We shall denote by $T^w$ the subgroup of $T$ that is centralized by any representative $\w$ of $w$.

We shall make extensive use of Chevalley's commutator formula (\cite[Theorem 5.2.2]{carter1}): for $\alpha$ and $\beta$ linearly independent roots and $a,\,b\in k$ there are structure constants $c^{ij}_{\alpha\beta}$ in the prime field of $k$ such that
\begin{equation}\label{chev}
x_\alpha(a)x_\beta(b)=x_\beta(b)x_\alpha(a)\prod_{i,\,j>0}x_{i\alpha+j\beta}(c^{ij}_{\alpha,\beta}a^ib^j)
\end{equation}
where the product is taken over all $(i,j)$ such that $i\alpha+j\beta\in\Phi$ and in any order for which $i+j$ is increasing. Moreover, $c^{ij}_{\alpha,\beta}\in\{\pm1,\pm2,\pm3\}$ and $3$ occurs only if $\Phi$ has a component of type $G_2$, so $c^{ij}_{\alpha\beta}\neq0$.\\

Given an element $x\in G$ we shall denote by $\O_x$ the conjugacy
class of $x$ in $G$ and by $G_x$
(resp. $B_x$, resp. $T_x$) the centralizer of $x$ in $G$ (resp. $B$,
resp. $T$). For a conjugacy class $\O=\O_x$ we shall denote by $\V$ the set of $B$-orbits
into which $\O$ is parted.

\begin{definition}\label{sferica} Let $K$ be a connected algebraic
  group. A homogeneous $K$-space
is called spherical if it has a dense orbit for some Borel subgroup of $K$.
\end{definition}

It is well-known (\cite{Bri}, \cite{Vin} in characteristic $0$,
\cite{gross}, \cite{knop} in positive characteristic)
that $X$ is a spherical homogeneous $G$-space if and only if  the set of
$B$-orbits in $X$ is finite.

\section{$B$-orbits and Bruhat decomposition}\label{sect-g2}

Let $\V$ be the set of $B$-orbits in a conjugacy class $\O$ in $G$. Since $G=\bigcup_{w\in W}BwB$ there is a natural map
$\phi\colon \V\to W$ associating to
$v\in\V$ the element $w$ in the Weyl group of $G$ for which
$v\subset BwB$.

\smallskip

It is shown in \cite{gio} for $G$ simple that if $\O$ is spherical as a homogeneous space then
the image of $\phi$ consists of involutions. The same proof holds for $G$ reductive. This motivates the
following definition.

\begin{definition}\label{quasisferica} Let $G$ be a connected reductive algebraic
  group. A conjugacy class $\O$ in $G$ is called
  quasi-spherical if the image of $\phi$ consists of involutions.
\end{definition}

\begin{remark}\label{regular-springer}{\rm Regular conjugacy classes
 in simple algebraic groups of rank greater than $1$ cannot be
 quasi-spherical. Indeed, by \cite[Theorem 8.1]{regular} regular classes meet Bruhat cells corresponding to Coxeter elements.}
\end{remark}

\subsection{The case of $G_2$}\label{g2}

We aim at showing that every quasi-spherical conjugacy class is spherical.
In order to illustrate this result explicitely, we analyze
quasi-spherical conjugacy classes for $G$ simple of type $G_2$ by inspection, making use of the classification of unipotent conjugacy classes (see, for instance, \cite[Section 7.12]{Hu-cc}) and of the commutator formula \eqref{chev}. 

\begin{theorem}Let $G$ be simple of type $G_2$. Then every quasi-spherical
  conjugacy class of $G$ is spherical.
\end{theorem}
\pf Let $\alpha$ and $\beta$ denote the short and long simple roots, respectively, and let $\O$ be a conjugacy class in $G$. We first assume that $\O$ is unipotent so it is either of type $A_1$, $\tilde{A}_1$, subregular or regular (\cite[Section 7.18]{Hu-cc}). If $\O$
is of type $A_1$ or $\tilde{A}_1$ then $\O$ is spherical, hence
quasi-spherical by \cite[Proposition 6, Proposition 11]{ccc} which
hold for arbitrary $k$. Alternativley, one may use
\cite[Theorem 3.2]{pany2} and \cite[Theorem 4.14]{FR}. If $\O$ is regular it cannot be
quasi-spherical by Remark \ref{regular-springer}.\\
The element
$u=x_{\beta}(1)x_{3\alpha+\beta}(1)\in G$ does not lie in the regular
unipotent conjugacy class by \cite[Lemma 3.2(c)]{regular}. Its conjugate
$\w_0u\w_0^{-1}=x_{-\beta}(a)x_{-3\alpha-\beta}(b)$ with $ab\neq0$ lies in $Bs_\beta
Bs_{3\alpha+\beta}B=Bs_\beta s_{3\alpha+\beta}B$ by \cite[Lemma 8.1.4(i), Lemma 8.3.7]{springer} and $s_\beta
s_{3\alpha+\beta}$ is not an involution. Then its class is not
quasi-spherical and by exclusion it is the subregular unipotent
conjugacy class, so the statement holds for unipotent conjugacy classes.

Let us now consider a representative $x\in\O\cap B$ with Jordan
decomposition $x=su\in TU$ with $s\neq1$. Then $G_s$ is connected and reductive (\cite[Theorem 2.2, Theorem 2.11]{Hu-cc}). We shall analyze the different cases according to the semisimple rank ${\rm srk}$ of $G_s$. 

If ${\rm srk}G_s=0$ then $\O$ is regular, hence it is not quasi-spherical by Remark \ref{regular-springer}.

If ${\rm srk}G_s=1$ and $u\neq1$ then $\O$ is regular,
hence it is not quasi-spherical. Let us assume $u=1$.
Up to conjugation by an
element in $N(T)$ we may assume that $G_s$ is
either $H_1=\langle T, X_{\pm \beta}\rangle$  or $H_2=\langle T, X_{\pm
  \alpha}\rangle$.

If $G_s=H_1$ conjugation of $s$ by
$x_{-\alpha}(1)x_{-\alpha-\beta}(1)$ yields
$$s_1=sx_{-\alpha}(a)x_{-\alpha-\beta}(b)x_{-2\alpha-\beta}(c)x_{-3\alpha-\beta}(d)x_{-3\alpha-2\beta}(e)\in
\O$$ for $a,b,c,d,e\in k$ with $ab\neq0$. Conjugation by a suitable element in $X_{-2\alpha-\beta}$ gives
$$s_2=sx_{-\alpha}(a)x_{-\alpha-\beta}(b)x_{-3\alpha-\beta}(d')x_{-3\alpha-2\beta}(e')\in
\O$$ for $d',e'\in k$. Conjugation by a suitable element in $X_{-3\alpha-\beta}$ gives
$$s_3=sx_{-\alpha}(a)x_{-\alpha-\beta}(b)x_{-3\alpha-2\beta}(e')\in
\O$$ and conjugation by a suitable element in $X_{-3\alpha-2\beta}$
gives
$$s_4=sx_{-\alpha}(a)x_{-\alpha-\beta}(b)\in
\O\cap B s_\alpha s_{\alpha+\beta} B$$ so $\O$ is not quasi-spherical.

Let $G_s=H_2$. Conjugation of $s$ by $x_{-3\alpha-\beta}(1)x_{-\beta}(1)$
yields
$$s_1=sx_{-\beta}(a)x_{-3\alpha-\beta}(b)x_{-3\alpha-2\beta}(c)\in \O$$
for some $a,b,c\in k$ with $ab\neq0$. Conjugation by a suitable element in
$X_{-3\alpha-2\beta}$ gives
$$s_2=sx_{-\beta}(a)x_{-3\alpha-\beta}(b)\in \O\cap B s_\beta
s_{3\alpha+\beta}B$$ so $\O$ is not quasi-spherical, concluding the analysis if ${\rm srk}G_s=1$.

Let ${\rm srk}G_s=2$ with $s\neq 1$. Up to conjugation by an element in $N(T)$ we may assume that $G_s$ is either
$$H_3=\langle T, X_{\pm \beta}, X_{\pm (3\alpha+\beta)}, X_{\pm
  (3\alpha+2\beta)}\rangle\mbox{ or }
H_4=\langle T, X_{\pm \beta}, X_{\pm
  (2\alpha+\beta)}\rangle.$$
If $G_s=H_3$ of type $A_2$ and $u=1$ then $\O$ is spherical by
\cite[Proposition 6, Theorem 16]{ccc} whose proofs hold in arbitrary good odd characteristic. Let us assume that $u\neq1$.
If $u$ is regular in $H_3$ then $\O$ is regular by \cite[Corollary 3.7]{regular}, hence it
is not quasi-spherical. It remains to analyze the class of $x=sx_{-\beta}(1)$. Conjugating by
$x_{-\alpha}(1)$ and reordering the terms gives
$$x_1=s
x_{-\beta}(1)x_{-\alpha-\beta}(b)x_{-2\alpha-\beta}(c)x_{-3\alpha-\beta}(d)x_{-3\alpha-2\beta}(e)x_{-\alpha}(f)\in
\O$$ for some $b,\,c,\,d,\,e,\,f\in k$ with $f\neq0$.
We can get rid of the term in $X_{-\alpha-\beta}$ conjugating by a
suitable element in $X_{-\alpha-\beta}$. Then we can get rid of the term in $X_{-2\alpha-\beta}$ conjugating by a
suitable element in $X_{-2\alpha-\beta}$ and, finally, we can get rid
of the term in $X_{-3\alpha-2\beta}$ by conjugating by a suitable
element in $X_{-3\alpha-\beta}$ obtaining
$$x_2=sx_{-\beta}(1)x_{-3\alpha-\beta}(b_1)x_{-\alpha}(f)\in\O$$ for some $b_1\in k$.
If $b_1=0$ then $\O\cap Bs_\beta s_\alpha B\neq\emptyset$ so $\O$ is
not quasi-spherical. If $b_1\neq0$ we have, for some $h\in T$ and some nonzero $a_i\in k$:
\begin{align*}
x_2&=shx_{\beta}(a_1)n_\beta
x_{\beta}(a_2)x_{3\alpha+\beta}(a_3)n_{3\alpha+\beta}x_{3\alpha+\beta}(a_4)x_{-\alpha}(f)\\
&=sh x_\beta(a_1)n_\beta x_{3\alpha+\beta}(a_3)x_\beta(a_2)x_{3\alpha+2\beta}(a_5)n_{3\alpha+\beta}x_{3\alpha+\beta}(a_4)x_{-\alpha}(f)\\
&\in TX_\beta X_{3\alpha+2\beta}n_\beta n_{3\alpha+\beta}X_{3\alpha+2\beta}X_\beta x_{-\alpha}(f)\\
&\subset Bn_\beta n_{3\alpha+\beta}P_\alpha^u x_{-\alpha}(f)\subset Bn_\beta n_{3\alpha+\beta}x_{-\alpha}(f)U\subset BX_{2\alpha+\beta} n_\beta n_{3\alpha+\beta} U
\end{align*}
so $\O\cap Bs_\beta s_{3\alpha+\beta} B\neq\emptyset$ and $\O$ is not
quasi-spherical.\\

Let $G_s=H_4$ be of type $A_1\times \tilde{A_1}$. If $u=1$ then
$\O$ is spherical by the argument in \cite[Theorem 16]{ccc}. If $u$ has
nontrivial components both in $A_1$ and in $\tilde{A}_1$ then $\O$ is
regular, hence it is not quasi-spherical. We are left with the analysis
of the classes of $y=sx_{-\beta}(1)$ and $z=sx_{-2\alpha-\beta}(1)$. \\
Conjugating $y$ by $x_{-3\alpha-\beta}(1)$ we get
$y_1=sx_{-3\alpha-\beta}(a)x_{-\beta}(1)x_{-3\alpha-2\beta}(b)$ for some $a,\,b\in k$ with
$a\neq0$. Conjugation by a suitable element in $X_{-3\alpha-2\beta}$
yields
$$y_2=sx_{-3\alpha-\beta}(a)x_{-\beta}(1)\in\O\cap
Bs_{3\alpha+\beta}s_\beta B$$ hence $\O_y$ is not quasi-spherical.

Conjugating $z$ by $x_{-\alpha}(1)$ we get
$z_1=sx_{-\alpha}(a)x_{-2\alpha-\beta}(1)x_{-3\alpha-\beta}(c)$ for some $a,\,c\in k$ with
$a\neq0$. Then conjugating $z_1$ by a suitable element in $X_{-\alpha-\beta}$
we obtain the element
$z_2=sx_{-\alpha}(a)x_{-\alpha-\beta}(d)x_{-3\alpha-\beta}(c_1)x_{-3\alpha-2\beta}(c_2)$ 
for some $c_1, c_2, d\in k$ with $d\neq0$. We can get rid of the term in $X_{-3\alpha-\beta}$
conjugating by a suitable element in $X_{-2\alpha-\beta}$ and then we can get rid of the term in $X_{-3\alpha-2\beta}$
conjugating by a suitable element in $X_{-3\alpha-2\beta}$.

Thus $z_3=sx_{-\alpha}(1)x_{-\alpha-\beta}(d)\in \O\cap Bs_\alpha
s_{\alpha+\beta}B$ hence $\O_z$ is not quasi-spherical. This exhausts
the list of conjugacy classes for $G$ of type $G_2$ and we have verified that all quasi-spherical conjugacy classes are spherical.\hfill$\Box$.

\section{Maximal $B$-orbits}\label{maximal}

Let $\O$ be a conjugacy class of $G$. Since $\O$ is an irreducible variety there exists a unique element in $W$ for which $\O\cap B wB$ is dense
in $\O$. We shall denote this element by $z_{\O}$. Denoting by $\overline{X}^Y$ the Zarisky closure of $X$ in $Y$ we have  
$$\O\subset \overline{\O}^G=\overline{\O\cap B z_{\O} B}^G\subset
\overline{B z_\O B}^G=\bigcup_{\sigma\leq z_\O}B\sigma B$$
so the element $z_\O$ is maximal in the image of $\phi$
 (cfr. \cite[Section 1]{ccc}).
We will call {\em maximal orbits} the elements $v$ in $\V$ for which
$\phi(v)=z_\O$ and we shall denote by $\V_{\max}$ the set of maximal
$B$-orbits in $\O$.

\begin{lemma}\label{finite}The following are equivalent for a conjugacy class $\O$ in $G$.
\begin{enumerate}
\item $\O$ is spherical.
\item $\V_{\max}$ contains only one element.
\item $\V_{\max}$ is a finite set.
 \end{enumerate}
\end{lemma}
\proof It follows from \cite[Corollary 26]{ccc}, \cite[Corollary 4.11]{gio} that if $\O$ is
spherical then $\V_{\max}$ contains only one element, namely the dense
$B$-orbit so 1 implies 2 and 2 trivially implies 3. Let us show that 3
implies 1. Since $\cup_{v\in \V_{\max}}v=\O\cap Bz_\O B$
is dense in $\O$ we have $\O\subset \cup_{v\in\V_{\max}}\overline{v}^{\O}$
with $\O$ irreducible (\cite[Proposition 1.5]{Hu-cc}) and $\V_{\max}$ a finite set. Then there necessarily exists $v_0\in \V_{\max}$ which is
dense in $\O$. \hfill$\Box$

\smallskip

Let us  analyze the maximal $B$-orbits in quasi-spherical
conjugacy classes.

\begin{lemma}\label{cnes-maximal}Let $\O$ be a quasi-spherical conjugacy class
with $w=z_\O$. Let $v\in \V_{\max}$ and let $x=u\w\v\in v$ with
  $u\in U^w$, $\w\in N(T)$ and $\v\in U$. Then for every
$\alpha\in\Delta$ such that $ws_\alpha>w$ in the Bruhat order we have:
\begin{enumerate}
\item  $s_{\alpha}w=ws_{\alpha}$ so $w\alpha=\alpha$;
\item $\v\in P^u_\alpha$, the unipotent radical of $P_\alpha$;
\item $X_{\pm \alpha}$ commutes with $\w$.
\end{enumerate}
\end{lemma}
\proof This is proved as  \cite[Lemma
  3.4]{gio}, since the proof therein uses only maximality of $w$ and
that $\O$ is quasi-spherical. \hfill$\Box$

\smallskip

\begin{lemma}Let $\O$ be a quasi-spherical conjugacy class with
  $w=z_\O$, let
  $\Pi=\{\alpha\in\Delta~|~w(\alpha)=\alpha\}$ and let $w_\Pi$ be the longest element in
  $W_\Pi$. Then $w=w_\Pi w_0$.
\end{lemma}
\proof By Lemma \ref{cnes-maximal}
if $\alpha\in\Delta$ and  $w\alpha\in\Phi^+$ then
$w\alpha=\alpha$. The statement follows
from  \cite[Proposition  3.5]{results}.\hfill$\Box$

\smallskip

The Lemmas above show that maximal $B$-orbits in quasi-spherical conjugacy classes behave similarly to
the dense $B$-orbit $v_0$ in a spherical conjugacy class. The analysis of
$z_\O$ given in \cite{gio} applies. 
\begin{proposition} The following properties hold for a
  quasi-spherical conjugacy class $\O$ with $w=z_\O=w_0w_\Pi$.
\begin{enumerate}
\item $\Pi$ is invariant with respect to $\vartheta=-w_0$;
\item The restriction of
$w_0$ to $\Phi(\Pi)$ coincides with $w_\Pi$;
\item $\Phi_w=\Phi\setminus\Phi(\Pi)$, notation as in \eqref{fivu};
\item $U_w=\langle X_\gamma~|~\gamma\in\Phi(\Pi)\cap\Phi^+\rangle$ and it normalizes $U^w$, notation as in \eqref{uvu};
\item $U_w$ commutes with $\w$ if $x=u\w\v\in\O\cap U^w N(T) U$.
\end{enumerate}
\end{proposition}
\pf The proof is as in \cite[Section 3]{gio}.\hfill$\Box$

\smallskip

In \cite{gio} an analysis of the possible $\Pi$ for
which $\phi(v_0)=w_0w_\Pi=z_\O$ for the dense $B$-orbit $v_0$ of a spherical
conjugacy class in a simple algebraic group was given. The proof of \cite[Lemma 4.1]{gio} can be
adapted to the case of maximal $B$-orbits in quasi-spherical conjugacy
classes, yielding the following statement.

\begin{lemma}\label{quali-no}Let $\O$ be a quasi-spherical conjugacy
  class  and let $w=w_0w_\Pi=z_\O$.
Let $\alpha$ and $\beta$ be simple roots with the following
properties:
$(\beta,\,\beta)=(\alpha,\,\alpha)$;
$w_0(\beta)=-\beta$;
$\beta\not\perp\alpha$;
$\beta\perp\alpha'$ for every
  $\alpha'\in\Pi\setminus\{\alpha\}$.

Then $\{\alpha\}$ cannot be a connected component of $\Pi$.
In particular, the list of the possible subsets
$\Pi$ for which $z_\O=w_0w_\Pi$ for $G$ simple coincides with the list
given in \cite[Corollary 4.2]{gio}.
\end{lemma}
\pf The proof follows as in \cite[Lemma 4.1]{gio} since it only
uses maximality of $w$ and that $\O$ is quasi-spherical. There, the proof is given for $G$
simple but it holds for $G$ reductive, too.\hfill$\Box$

\smallskip

Let  $\O$ be
quasi-spherical with $w=z_\O=w_0w_\Pi$ and
let $\Phi_1=\Phi\cap {\rm Ker}(1+w)$. Then $\Phi_1$ is a root subsystem
of $\Phi$ and we put $\Phi_1^+=\Phi^+\cap\Phi_1$. If we write $w=\prod_js_{\gamma_j}$ as a
product of reflections with respect to mutually orthogonal roots
then each $\gamma_j$ lies in $\Phi_1$. We shall denote by $W(\Phi_1)$ the subgroup of $W$ generated by reflections with respect to roots in $\Phi_1$, so $w\in W(\Phi_1)$.

\begin{lemma}Let notation be as above and let $\beta\in\Phi$. Then $\beta\in\Phi_1$ if and only if $\beta\perp \Pi$ and $\vartheta\beta=\beta$.
\end{lemma}
\pf  If $\beta\perp\Pi$ and $\vartheta\beta=\beta$ then
$w_\Pi\beta=\beta$ and $w_0\beta=-\beta$ thus $\beta\in\Phi_1$.\\ 
Conversely, if $w\beta=-\beta$ then for every $\alpha\in\Pi$ we have
$\beta\perp\alpha$ because $\alpha$ and $\beta$ lie in distinct
eigenspaces of the orthogonal transformation $w$. 
Let now $\alpha\in\Phi$ and $w\beta=-\beta$. We have 
$$(\vartheta\beta,\alpha)=-(w_0\beta,\alpha)=-(w
w_\Pi\beta,\alpha)=-(w\beta,\alpha)=(\beta,\alpha)$$ and since this
holds for every $\alpha$, we have the statement.\hfill$\Box$

\smallskip

Let us denote by $G(\Phi_1)$ the subgroup of $G$ generated by
$T$ and the root subgroups $X_{\pm\beta}$ with $\beta\in\Phi_1$. Let
$U_{\Phi_1}=\langle X_\beta,\,\beta\in\Phi_1^+\rangle$.

The following Lemma is an analogue of \cite[Lemma 4.8, Remark
  4.9]{gio} for quasi-spherical conjugacy classes.

\begin{lemma}\label{ortogonale}Let $\O$ be a quasi-spherical conjugacy
  class and let $z_\O=w_0w_\Pi$. Let $\w\in N(T)$ be a representative of $w$.
Then for every $x=\w t\v \in \w B\cap \O$ we have $\v\in U_{\Phi_1}$, $w\in
W(\Phi_1)$ and $x$ commutes with  $(T^w)^\circ U_w$.
\end{lemma}
\pf The proof when $G$ is simple follows exactly as in \cite[Lemmas 4.5, 4.6, 4.7, 4.8,
  4.9]{gio}. 
Indeed, for their proofs 
we only need $w$ to be maximal,  $\O$ to be quasi-spherical and the list in \cite[Corollary
  4.2]{gio}. The general case is a consequence of the case of $G$ simple.\hfill$\Box$

\smallskip

\begin{lemma}\label{U-}Let $\O$ be a quasi-spherical conjugacy class
  and let $w=w_0w_\Pi=z_\O$. Then $\langle
  X_{-\alpha}, \alpha\in\Pi\rangle$ commutes with every $x=\w t\v \in
  \w B\cap \O$.
\end{lemma}
\Pf It is not restrictive to assume $G$ to be simple. By Lemmas \ref{cnes-maximal} and \ref{ortogonale} it is enough
to show that $X_{-\alpha}$ commutes with $\v$ for every
$\alpha\in\Pi$. If this were not the case, by \eqref{chev} there would occur in
the expression of $\v$ at
least one root subgroup $X_\gamma$ with nontrivial coefficient and
with $\gamma-\alpha\in\Phi$. We consider such a $\gamma$ of minimal height.
By Lemma \ref{ortogonale} and \cite[Chapitre 6, \S 1.3]{bourbaki}
this could happen only if $\Phi$ is doubly-laced and $\alpha$ is a
short root. Then we would also have
$\alpha+\gamma\in\Phi$, which is impossible because $X_\alpha$ 
commutes with $\v$ by Lemma \ref{ortogonale}. \hfill$\Box$ 

\smallskip

A consequence of Lemma \ref{ortogonale} is the following result.
\begin{proposition}\label{pro}Let $\O$ be a quasi-spherical conjugacy
  class, let $w=z_\O=w_0w_\Pi$ and let $v\in\V_{\max}$. Then $\dim(v)=\ell(w)+{\rm rk}(1-w).$
\end{proposition}
\pf Let $n$ be the rank of $G$ and let $x=\w\v\in v$. By Lemma
\ref{ortogonale} the centralizer $B_x$ of
$x$ in $B$ contains $(T^w)^\circ U_w$. On the other hand, if
$b=u^wu_wt\in U^wU_wT$ commutes with $x$ we have
$$
\w\v u^wu_wt =u^wu_wt \w \v=u^w u_w \w (\w^{-1}t\w)\v=u^w\w u_w
(\w^{-1}t\w)\v
$$
where for the last equality we used Lemma \ref{cnes-maximal}. By
uniqueness in the Bruhat decomposition we have $u^w=1$ so $B_x\subset
T_xU_w$ because $U_w\subset B_x$. Moreover, if $t\in T_x$ we have
$$
\w (\w^{-1}t\w)\v=t\w\v=\w\v t \in \w TU
$$
and uniqueness of the decomposition in $TU$ gives $t\in
T^w$. Therefore $(T^w)^\circ U_w \subset B_x\subset T^wU_w$ and
$\dim v=|\Phi^+|+n-(|\Phi^+|-\ell(w))-(n-{\rm rk}(1-w))=\ell(w)+{\rm
  rk}(1-w)$. \hfill$\Box$

\section{The case $z_\O=w_0=-1$}\label{bigcell}

In this section $\Phi$ is such that $w_0$ acts as $-1$ in the geometric
representation of $W$. If $\O$ is a quasi-spherical conjugacy class
intersecting the big Bruhat cell $Bw_0B$ then $\Phi_1=\Phi$ and
$\Pi=\emptyset$ so Lemma \ref{ortogonale} gives no restriction to a representative $x=\w \v\in
\O\cap\w U$. For this reason we use a different approach for such
classes. 
%
%
%
By Lemma \ref{finite} if a conjugacy class has finitely-many maximal $B$-orbits then it is spherical. The aim of this Section is to show that every quasi-spherical conjugacy class $\O$ intersecting $B w_0 B$ has only finitely-many maximal $B$-orbits. This will be achieved by counting the possible representatives of a maximal $B$-orbit lying in $\w_0 U$ for a fixed $\w_0\in N(T)$. Next Lemma shows that every maximal $B$-orbit meets $\w_0 U$.      

\begin{lemma}\label{meet}Let $G$ be simple and let $\O$ be a quasi-spherical conjugacy class with
  $z_\O=w_0=-1$. For any $v\in \V_{\max}$ and any representative $\w_0$
  of $w_0$ in $N(T)$ we have $v\cap \w_0 U\neq\emptyset$.
\end{lemma}
\pf Let $x=u\w_0 t\v\in v\cap U\w_0 B$. Then for every $s\in T$ we have
$x_s=s^{-1}u^{-1}u\w_0 t\v us=\w_0 s^2t u'\in v\cap \w_0 TU$ and since the map
$s\mapsto s^2\in T$ is onto (\cite[III.8.9]{borel}) we may choose $s$ so that
$x_s\in v\cap\w_0 U$.\hfill$\Box$

\smallskip

\begin{lemma}\label{alpha-beta}Let $\O$ be a quasi-spherical conjugacy class with  $z_\O=w_0=-1$. Let $\w_0$ be a representative of $w_0$ and let $x=\w_0\v\in\O\cap\w_0 U$, with
  $\v=\prod_{\gamma\in\Phi^+}x_\gamma(c_\gamma)$ in a fixed ordering. Let $\alpha$ and $\beta$ be adjacent simple
  roots of the same length. Then the number of possibilities for
  $c_\alpha$ and $c_\beta$ is finite and $c_{\alpha+\beta}$ is
  completely determined by the ordering, $c_\alpha$ and $c_\beta$.
\end{lemma}
\pf Let $P=P_{\{\alpha,\beta\}}$ with unipotent radical $P^u$. Let us assume that $\alpha$ precedes $\beta$ in the
ordering. We may write: $x=\w_0\v\in\w_0
x_{\alpha}(c_\alpha)x_{\beta}(c_\beta)x_{\alpha+\beta}(c_{\alpha+\beta})P^u$.
For $h\in k$ we put
$y(h):=n_\alpha x_{\alpha}(h)x x_{\alpha}(-h)n_\alpha^{-1}$. Then, for
some nonzero structure constants $\theta_1,\theta_2,\theta_3, c^{11}_{\alpha\beta}$ and some $t_1\in T$ we have
\begin{align*}
y(h)&\in n_\alpha\w_0 x_{-\alpha}(\theta_1
h)x_{\alpha}(c_\alpha-h)x_\alpha(h)x_{\beta}(c_\beta)x_{\alpha+\beta}(c_{\alpha+\beta})x_{\alpha}(-h)n_\alpha^{-1}P^u\\
&=\w_0 t_1 x_\alpha(\theta_1\theta_2
h)x_{-\alpha}(\theta_3(c_\alpha-h))n_\alpha
x_{\beta}(c_\beta)x_{\alpha+\beta}(c_{\alpha+\beta}+hc_\beta
c^{11}_{\alpha\beta})n_\alpha^{-1}P^u.
\end{align*}
Let $h_1$ and $h_2$ be the solutions of
\begin{equation*}X^2(\theta_1\theta_2\theta_3)-c_\alpha
\theta_1\theta_2\theta_3X-1=0\end{equation*} so that $-(\theta_1\theta_2
h_i)^{-1}=(c_\alpha-h_i)\theta_3$. By \cite[Lemma 8.1.4 (i)]{springer} we have
\begin{align*}
y(h_i)&\in\w_0 t_1 n_\alpha t_2
x_{\beta+\alpha}(c'_\beta)x_{\beta}(\theta_4(c_{\alpha+\beta}+h_ic_\beta
c^{11}_{\alpha\beta}))P^u\\
&\subset \w_0n_\alpha t_3 x_{\beta}(\theta_4(c_{\alpha+\beta}+h_ic_\beta
c^{11}_{\alpha\beta}))P_\beta^u\subset \O\cap Bw_0s_\alpha B
\end{align*}
for some $t_2, t_3\in T$, some $c_\beta'\in k$ and some nonzero structure constant
$\theta_4$.
Since $w_0s_{\alpha+\beta}\beta=\alpha\in\Phi^+$, conjugation of $y(h_i)$ by $n_\beta$ would yield an element in $\O\cap
Bw_0s_{\alpha+\beta} s_\beta B$ unless
\begin{equation}\label{condiz}c_{\alpha+\beta}+h_ic_\beta
c^{11}_{\alpha\beta}=0.\end{equation} As
$s_{\alpha+\beta}s_\beta$ is not an involution, \eqref{condiz} must hold for both $i=1,2$ thus we have either
$h_1=h_2$ so that
\begin{eqnarray}\label{Ialpha}&\Delta_\alpha=\theta^2_1\theta^2_2\theta^2_3c_\alpha^2+4\theta_1\theta_2\theta_3=0,\quad\mbox{
or }\\
\label{IIalpha}&c_\beta=c_{\alpha+\beta}=0.\end{eqnarray}

Let us now consider, for $l\in k$, the element
\begin{align*}
z(l)&=n_\beta x_\beta(l)x x_{\beta}(-l)n_\beta^{-1}\\
&\in n_\beta
x_\beta(l)\w_0
x_\beta(c_\beta)x_\alpha(c_\alpha)x_{\alpha+\beta}(c_{\alpha+\beta}+
c_\alpha c_\beta c^{11}_{\alpha\beta}) x_{\beta}(-l)n_\beta^{-1}P^u.
\end{align*}
Repeating the same argument for $\beta$ we see that there are nonzero
structure constants $\eta_1,\eta_2,\eta_3,\eta_4$ so that if
$l_j$ is a solution of
\begin{equation*}\eta_1\eta_2\eta_3 X^2-c_\beta
\eta_1\eta_2\eta_3 X-1=0\end{equation*} then
$$z(l_j)\in \O\cap\w_0 n_\beta T x_\alpha(\eta_4(c_{\alpha+\beta}+c_\alpha
c_\beta c_{\alpha\beta}^{11}-l_jc_\alpha
c_{\alpha\beta}^{11}))P_\alpha^u$$
so conjugation by $n_\alpha$ would yield an element in $\O\cap Bw_0
s_{\alpha+\beta} s_\alpha B$ unless
\begin{equation}\label{alpha+beta}c_{\alpha+\beta}+c_\alpha
c_\beta c_{\alpha\beta}^{11}-l_jc_\alpha
c_{\alpha\beta}^{11}=0\end{equation}
for both $j=1,2$. This forces either $l_1=l_2$ and therefore
\begin{eqnarray}\label{Ibeta}&\Delta_\beta=\eta^2_1\eta^2_2\eta^2_3c_\beta^2+4\eta_1\theta_2\theta_3=0,\quad\mbox{
    or }\\
\label{IIbeta}&c_\alpha=c_{\alpha+\beta}=0.
\end{eqnarray}
If \eqref{Ialpha} holds then $c_\alpha\neq0$ so \eqref{Ibeta} holds.
Then the possibilities for $c_\alpha$ and $c_\beta$ are finite. Besides, by
\eqref{alpha+beta} the coefficient
$c_{\alpha+\beta}$ is completely determined by $c_\alpha$, $c_\beta$, and the structure
constants.\\
If \eqref{IIalpha} holds then \eqref{Ibeta} cannot hold so
$c_\alpha=c_\beta=c_{\alpha+\beta}=0$, whence the statement.\hfill$\Box$

\smallskip

\begin{lemma}\label{a<=b}Let $\O$, $z_\O$, $x$, $\v$ be as in Lemma \ref{alpha-beta}. Let $\alpha$ and $\beta$ be adjacent simple
  roots with $2(\alpha,\alpha)=(\beta,\beta)$. Then
  $c_{2\alpha+\beta}$ and $c_{\alpha+\beta}$ are completely determined
  by the ordering, $c_\alpha$ and $c_\beta$.
\end{lemma}
 \pf It is not restrictive to assume that $\alpha$ precedes $\beta$ and $\alpha+\beta$ in
 the ordering. Let $P=P_{\{\alpha,\beta\}}$ and $P^u$ be its unipotent
 radical. Then $$x\in\w_0
 x_{\alpha}(c_{\alpha})x_{\beta}(c_{\beta})x_{\alpha+\beta}(c_{\alpha+\beta})x_{2\alpha+\beta}(c_{2\alpha+\beta})P^u.$$
 Conjugation by $n_\alpha x_\alpha(h)$ for $h\in k$ yields 
$$y(h)\in\w_0t_1 x_{\alpha}(\eta_1
 h)x_{-\alpha}(\eta_2(c_\alpha-h))x_{\beta}(\eta_3(c_{2\alpha+\beta}+hc_{\alpha+\beta}c^{11}_{\alpha,\alpha+\beta}+c_\beta h^2c_{\alpha,\beta}^{21}))P_\beta^u$$
for some $t_1\in T$ and some nonzero structure constants
 $\eta_1,\eta_2,\eta_3$.
If $h_1,\,h_2$ are the solutions of
\begin{equation*}\eta_1\eta_2X^2-c_\alpha
\eta_1\eta_2X-1=0
\end{equation*}
then $y(h_1)$, $y(h_2)$ lie in $\O\cap Bw_0s_\alpha B$ and $n_\beta y(h_i)
n_\beta^{-1}\in\O\cap Bw_0s_{\alpha+\beta}s_\beta B$ unless
\begin{equation}\label{2alfa+beta}c_{2\alpha+\beta}+h_ic_{\alpha+\beta}c^{11}_{\alpha,\alpha+\beta}+c_\beta h_i^2c_{\alpha,\beta}^{21}=0
\end{equation}
for both $i=1,2$. Besides, $h_1+h_2=c_\alpha$ and $h_1h_2=-(\eta_1\eta_2)^{-1}$. Thus we have either
\begin{eqnarray}\label{1alpha}&\Delta_\alpha=(\eta_1\eta_2c_\alpha)^2+4\eta_1\eta_2=0\mbox{
    and }
    c_{2\alpha+\beta}=ac_\alpha
    c_{\alpha+\beta}+bc_\beta c_\alpha^2
    \\
\label{2alpha}&\mbox{or } c_{2\alpha+\beta}=cc_\beta
   \mbox{ and }
c_{\alpha+\beta}=d c_\beta c_\alpha
\end{eqnarray}
for $a$, $b$, $c$, $d\in k$. On the other hand, reordering terms we have:
$$x\in\w_0
x_\beta(c_\beta)x_\alpha(c_{\alpha})x_{\alpha+\beta}(c_{\alpha+\beta}+c_\alpha
c_\beta
  c^{11}_{\alpha\beta})x_{2\alpha+\beta}(c_{2\alpha+\beta}+c_\alpha^2c_\beta c_{\alpha\beta}^{21})P^u.$$
Conjugation by $n_\beta x_\beta(l)$ for $l\in k$ gives an element
$$z(l)\in\w_0 t_2 x_\beta(\theta_1
l)x_{-\beta}(\theta_2(c_\beta-l))x_{\alpha}(\theta_3(c_{\alpha+\beta}+c_\alpha
c_\beta
c_{\alpha\beta}^{11}-lc_\alpha c_{\alpha\beta}^{11}))P_\alpha^u$$
for some $t_2\in T$ and some nonzero structure constants
$\theta_1,\theta_2,\,\theta_3,\,c_{\alpha\beta}^{11}$. For the solutions
$l_1,\,l_2$ of
\begin{equation*}\theta_1\theta_2X^2-c_\beta
\theta_1\theta_2X-1=0
\end{equation*}
the corresponding elements $z(l_1)$, $z(l_2)$ lie in $\O\cap Bw_0s_\beta B$ and $n_\alpha z(l_i)
n_\alpha^{-1}\in Bw_0s_{2\alpha+\beta}s_\alpha B$ unless
\begin{equation}\label{alfa+beta}c_{\alpha+\beta}+c_\alpha
c_\beta c_{\alpha\beta}^{11}-l_ic_\alpha c_{\alpha\beta}^{11}=0
\end{equation}
for $i=1,2$. It follows that we have either
\begin{eqnarray}\label{1beta}&\Delta_\beta=(\theta_1\theta_2c_\beta)^2+4\theta_1\theta_2=0\mbox{
    and }
    c_{\alpha+\beta}=d c_\alpha c_\beta\\
\label{2beta}&\mbox{ or } c_{\alpha+\beta}=c_\alpha=0.
\end{eqnarray}
Arguing as in Lemma \ref{alpha-beta} we see that $c_{\alpha+\beta}$ and $c_{2\alpha+\beta}$ are completely
determined by the ordering, $c_\alpha$ and $c_\beta$. \hfill$\Box$

\smallskip

\begin{lemma}\label{determined}Let $\Phi$ be a simply- or doubly-laced
  root system for which $w_0=-1$. Let $\O$, $z_\O$, $x$, $\v$ be as in Lemma \ref{alpha-beta}.
Then, for every $\gamma=\sum_{j\in J} n_j\alpha_j\in\Phi^+$ with $J\subset\{1,\ldots,\,n\}$  there is a polynomial $p_\gamma(X)\in k[x_j,\,j\in J]$ 
depending only on the fixed ordering of the positive roots and the structure constants of $G$ 
such that the coefficient $c_\gamma$ in the expression of $\v$ is the evaluation of $p_\gamma(X)$ at $x_j=c_{\alpha_j}$ for every $j\in J$. 
\end{lemma}
\pf We shall proceed by induction on the height ${\rm ht}$ of the root $\gamma$.
Let us assume that the claim holds for all $\gamma$ with ${\rm ht}(\gamma)\leq m-1$.
By Lemmas \ref{alpha-beta} and \ref{a<=b} the statement holds for $m=1,2$ so we assume $m\geq 3$.

Let $\nu\in\Phi^+$ with ${\rm
  ht}(\nu)=m$. Then there exists $\beta\in\Delta$ for which
${\rm ht}(s_\beta\nu)\leq m-1$. We put
\begin{equation}y=n_\beta x n_\beta^{-1}=\w_0 t
\prod_{\gamma\in\Phi^+}x_{s_\beta\gamma}(\theta_\gamma c_\gamma)\end{equation}
 for some nonzero structure constants $\theta_\gamma$. Here the products
  have to be intended in the fixed ordering of the $\gamma$'s. We have:
$$y=\w_0 t (\prod_{\gamma<_o\beta}x_{s_\beta\gamma}(\theta_\gamma
c_\gamma))x_{-\beta}(\theta_\beta c_\beta) (\prod_{\gamma>_o\beta}x_{s_\beta\gamma}(\theta_\gamma
c_\gamma))$$ where $<_o$
indicates that a root precedes another in the fixed ordering and the expression
makes sense also if $c_\beta=0$. Then
\begin{align*}
y&=\w_0 t x_{-\beta}(\theta_\beta c_\beta)x_{-\beta}(-\theta_\beta c_\beta)(\prod_{\gamma<_o\beta}x_{s_\beta\gamma}(\theta_\gamma
c_\gamma))x_{-\beta}(\theta_\beta c_\beta) (\prod_{\gamma>_o\beta}x_{s_\beta\gamma}(\theta_\gamma
c_\gamma))\\
&=x_\beta(\eta c_\beta)\w_0t (x_{-\beta}(-\theta_\beta c_\beta)(\prod_{\gamma<_o\beta}x_{s_\beta\gamma}(\theta_\gamma
c_\gamma))x_{-\beta}(\theta_\beta c_\beta))(\prod_{\gamma>_o\beta}x_{s_\beta\gamma}(\theta_\gamma
c_\gamma))
\end{align*}
for some nonzero structure constant $\eta$. Conjugation by
$x_{\beta}(-\eta c_\beta)$ yields 
\begin{align*}
z&=\w_0t (x_{-\beta}(-\theta_\beta c_\beta)(\prod_{\gamma<_o\beta}x_{s_\beta\gamma}(\theta_\gamma
c_\gamma))x_{-\beta}(\theta_\beta c_\beta))(\prod_{\gamma>_o\beta}x_{s_\beta\gamma}(\theta_\gamma
c_\gamma))x_\beta(\eta c_\beta)\\
&=\w_0t \prod_{\gamma\in\Phi^+}x_\gamma(d_\gamma)\in\w_0 tU\cap\O
\end{align*}
where the last equality indicates reordering of root subgroups. 
By the induction hypothesis applied to $z$ and $s_\beta\nu$, the coefficient $d_{s_\beta\nu}$ is evaluation at the 
 $d_\alpha$ for $\alpha$ in the support of $s_\beta\nu$  of a polynomial $q_{s_\beta\nu}(X)$.
Besides, each
$d_{\mu}$ differs from $\theta_{s_\beta\mu} c_{s_\beta\mu}$ by a
 (possibly trivial) sum of monomials in the $\theta_{\mu'}
 c_{\mu'}$, $c_\beta$ and the structure constants coming from application of
 \eqref{chev} when reordering root subgroups. More
 precisely, we have 
\begin{equation}\label{reordering}d_\mu=\theta_{s_\beta\mu} c_{s_\beta\mu}+\sum *(\prod_{l=1}^pc_{\nu_l}^{i_l})c_\beta^j
\end{equation}
where $*$ denotes a coefficient depending on the structure constants and the sum is taken over the possible decompositions 
$\mu=\sum_{l=1}^pi_ls_\beta\nu_l+j\beta$ for $i_l>0$ and $j\geq0$. 
In particular, if $\mu$ is simple there is no such decomposition: in this case
 $d_\mu=\theta_{s_\beta\mu} c_{s_\beta\mu}$ and by Lemmas \ref{alpha-beta} and \ref{a<=b} the coefficient $c_{s_\beta\mu}$ is evaluation of a polynomial at the $c_\alpha$ for $\alpha$ in the support of $s_\beta\mu$, and such support is contained in $\{\beta,\,\mu\}$. Thus, by the induction hypothesis 
 $d_{s_\beta\nu}$ is evaluation of a polynomial $\overline{q}(X)$ at the $c_\alpha$ for
 $\alpha$ in the support of $\nu$. We wish to prove that
 the same holds for $c_\nu$. Contribution to $d_{s_\beta\nu}$ as in \eqref{reordering} may occur when
 \begin{equation}\label{nu}s_\beta\nu=\sum_{l=1}^pi_ls_\beta\nu_l+j\beta\end{equation} for $i_l>0$ and
 $j\geq0$. Then ${\rm ht}(s_\beta\nu_l)<{\rm ht}(s_\beta\nu)\leq m-1$.
We wish to show that ${\rm ht}(\nu_l)\leq m-1$ for every $l$ so we may
 apply the induction hypothesis to $c_{\nu_l}$. Suppose that there is a decomposition \eqref{nu} and an $l$ for which ${\rm ht}(\nu_l)\geq m$ for some $l$. Since 
$$m\leq {\rm ht}(\nu_l)\leq {\rm ht}(s_\beta\nu_l)+2\leq {\rm ht}(s_\beta\nu)-1+2\leq m$$
we would necessarily have ${\rm ht}(\nu)={\rm ht}(\nu_l)=m$; ${\rm ht}(s_\beta\nu_l)=m-2$; ${\rm ht}(s_\beta \nu)=m-1$ thus $s_\beta\nu=s_\beta\nu_l+\alpha$ for some $\alpha\in\Delta$. Applying $s_\beta$ to this equality we would have $\nu=\nu_l+s_\beta\alpha$ contradicting ${\rm ht}(\nu)={\rm ht}(\nu_l)$. 


Thus induction applies and 
$$c_\nu=*\overline{q}(c_\alpha)_{\alpha\in{\rm supp}(\nu)}+\sum *\prod_{l=1}^p(p_{\nu_l}(c_\alpha)_{\alpha\in{\rm supp}(\nu_l)})^{i_l}c_\beta^j$$
is evaluation of a polynomial depending only on the structure constants.
\hfill$\Box$
\smallskip

\begin{remark}{\rm The proof of Lemma \ref{determined} can be adapted
  to show that if $\Phi$ is simply-laced and $\v\in P_\alpha^u$ for
  some $\alpha\in\Delta$ then $\v=1$ so $x=\w_0$ and $\O$ is a
  symmetric space.}
\end{remark}

\begin{proposition}\label{bigcellDE}Let $\Phi$ be irreducible, simply-laced, with $w_0=-1$ and let $\O$ be a
  quasi-spherical conjugacy class with $z_\O=w_0$. Then $\O$ is
  spherical.
\end{proposition}
\pf If $\Phi=\{\pm\alpha\}$ is of type $A_1$ the statement follows from even
dimensionality of conjugacy classes (\cite[Prop. 4.3]{even}). Indeed, $\dim\O\leq 2$ and for a representative $x$ of a maximal $B$-orbit in $\O$ we have $\dim
B.x=\ell(s_\alpha)+{\rm rk}(1-s_\alpha)=2$ by Proposition \ref{pro}. Thus, $B.x$ is dense in $\O$. 
 
We assume now that the rank of $\Phi$ is at least two.
By Lemma \ref{meet} every $v\in \V_{\max}$ meets $\w_0 U$ for
every choice of $\w_0$ and by Lemmas \ref{alpha-beta} and
\ref{determined} there is only a finite number of elements in $\O\cap
\w_0 U$. We conclude using Lemma \ref{finite}. \hfill$\Box$

\smallskip

\begin{proposition}\label{bigcellf4}Let $G$ be simple of type $F_4$ and let $\O$ be a
  quasi-spherical conjugacy class with $z_\O=w_0$. Then $\O$ is spherical.
 \end{proposition}
\pf By Lemma \ref{finite} we need to show that there are only finitely many maximal $B$-orbits. By Lemma \ref{meet} it is enough to show that there are only finitely many elements in $\w_0 U\cap\O$ for a fixed $\w_0\in N(T)$. By Lemma \ref{determined} it is enough to show that there are only finitely many possibilities for $c_{\alpha}$ for $\alpha\in\Delta$. Applying Lemma
\ref{alpha-beta} to the pair $\alpha_1,\,\alpha_2$
we see that for $x=\w_0\v\in\O\cap \w_0 U$ there is a finite number of
possibilities for the coefficients $c_{\alpha_1}$ and $c_{\alpha_2}$ in $\v$.  Applying Lemma
\ref{alpha-beta} to the pair $\alpha_3,\,\alpha_4$
we see that there is a finite number of
possibilities for $c_{\alpha_3}$ and $c_{\alpha_4}$, concluding the proof. 
\hfill$\Box$

\smallskip

\begin{proposition}\label{cn}Let $\Phi$ be irreducible of type $C_n$ and let $\O$
  be a quasi-spherical conjugacy class with $z_\O=w_0$. Then $\O$ is spherical.
\end{proposition}
\pf If $n=2$ by Proposition \ref{pro} we have $6=\dim B=\dim
v$ for every maximal $B$-orbit $v$. On the other hand $\dim\O< 2|\Phi^+|=8$ because $\O$ cannot be
regular (see Remark \ref{regular-springer}). It follows from even
dimensionality of conjugacy classes (\cite[Prop. 4.3]{even}) that $\dim \O=\dim v$ so $v$ is dense and $\O$ is spherical. 

Let us now assume that $n\geq3$. 
Let $\w_0\in N(T)$ be fixed and let
$x=\w_0\v=\w_0\prod_{\gamma\in\Phi^+}x_\gamma(c_\gamma)$ be as in Lemma \ref{alpha-beta}. By Lemmas \ref{finite}, \ref{meet} and \ref{determined}
it is enough to prove that there is a finite number of possibilities for $c_\alpha$ for
$\alpha\in\Delta$. By Lemma \ref{alpha-beta} we have either  $\Delta_{\alpha_i}=0$ for
$i=1,\ldots,n-1$ or
$c_{\alpha_i}=0$ for $i=1,\ldots,n-1$. In the first case, 
Lemma \ref{a<=b} with $\alpha=\alpha_{n-1}$
and $\beta=\alpha_n$ gives $\Delta_{\alpha_n}=0$ so there are finitely many possibilities for all $c_\gamma$. We shall thus focus on the case $c_{\alpha_i}=0$ for $i\leq n-1$. Then Lemma \ref{a<=b} with $\alpha=\alpha_{n-1}$
and $\beta=\alpha_n$ gives $c_{\alpha+\beta}=0$. We claim that $c_\gamma=0$ for
every short root. We proceed by induction as in Lemma \ref{determined} and we look at possible the contribution as in \eqref{reordering} to $c_\nu$ with $\nu=s_\beta\mu$ and ${\rm ht}(\mu)<{\rm ht}(\nu)$. This would correspond to a decomposition of the short root $s_\beta\nu=\sum i_js_\beta\nu_j+i\beta$ with $i_j>0$ and $i\geq0$.
If $i>0$ we have nontrivial contribution only if $\beta$ is a long root, for $c_\beta=0$ if $\beta$ is short. Thus, both for $i=0$ and $i>0$  there is at least one $\nu_j$ which is short and then ${\rm ht}(s_\beta\nu_j)\leq {\rm ht}(s_\beta\nu)-1=m-2$ so ${\rm ht}(\nu_j)\leq m-2+1$. By the induction hypothesis $c_{\nu_j}=0$ and there is no contribution coming from this decomposition, so the claim is proved. 

In other words, putting 
$\gamma_n=\alpha_n$ and
$\gamma_i=s_{i}s_{i+1}\cdots s_{n-1}\alpha_n$ for $i=1,\ldots,\,n-1$
we have $x=\w_0 \prod_{i=1}^nx_{\gamma_i}(a_i)$ for some
$a_i\in k$. We claim that there can be only finitely many elements of
this type in a fixed class $\O$. It is not restrictive to assume that
$G=Sp_{2n}(k)$. Then, $G$ is the subgroup of $GL_{2n}(k)$ of matrices preserving the bilinear form associated with the matrix $\begin{pmatrix}0&I\cr
-I&0\cr\end{pmatrix}$ with respect to the canonical basis of $k^{2n}$. We choose $B$ as the subgroup of $G$ of matrices of the form $\begin{pmatrix}X&XA\cr
0&^tX^{-1}\cr\end{pmatrix}$ where $X$ is an invertible upper triangular matrix, $^tX^{-1}$ is its inverse transpose and $A$ is a symmetric matrix. Then the computations above translate into:
$$x=x(D,A)=\w_0\v=\begin{pmatrix}0&D\cr
-D^{-1}&0\cr\end{pmatrix}\begin{pmatrix}I&A\cr 0&I\cr\end{pmatrix}=\begin{pmatrix}0&D\cr
-D^{-1}&-D^{-1}A\cr\end{pmatrix}$$ for some diagonal matrices
$D$ and $A$ with $D$ fixed and invertible.
It is immediate to verify that
for two diagonal matrices $A$ and $A'$, the characteristic polynomials
of $x(D,A)$ and  $x(D,A')$ are the same only if  $D^{-1}A$ and $D^{-1}A'$ coincide up to a permutation of the diagonal entries. Therefore there are only finitely
many matrices of the form $x(A,D)$ in a single conjugacy class $\O$.  Since by Lemma \ref{meet} each maximal $B$-orbit in $\O$ contains some $x(A,D)$ or some of the finitely many representatives with $\Delta_{\alpha_j}=0$ for every $j$, there are only finitely many maximal $B$-orbits in $\O$ and we may conclude by using Lemma \ref{finite}.
\hfill$\Box$

\begin{proposition}\label{bn}Let $\Phi$ be irreducible of type $B_n$ and let $\O$
  be a quasi-spherical conjugacy class with $z_\O=w_0$. Then $\O$ is spherical.
\end{proposition}
\pf  The case $n=2$ is dealt with in Proposition \ref{cn} so we may assume $n\geq3$.
Let $\w_0\in N(T)$ be fixed and let
$x=\w_0\v$ be as in Lemma \ref{alpha-beta}.
We shall show that there is a finite number of possibilities for
$c_\alpha$ for $\alpha\in\Delta$. It follows from Lemma
\ref{alpha-beta} that  we have either $c_\alpha=0$ for every long
simple root $\alpha$ or $\Delta_\alpha=0$ for every long
simple root $\alpha$. In the first case, Lemma \ref{a<=b} with 
$\alpha=\alpha_n$ and $\beta=\alpha_{n-1}$ shows that \eqref{1beta}
cannot be satisfied so $c_{\alpha_n}=0$ as well. Hence, there is no freedom for the $c_{\alpha}$ in this case and we shall focus on the second case. 
Let $\Delta'=\{\alpha_1,\,\ldots,\,\alpha_{n-1}\}$ and $P=P_{\Delta'}$. Then $x=\w_0\v_1 \v_2$ with
$\v_1\in \langle X_\alpha,\,\alpha\in\Delta'\rangle$ and $\v_2\in P^u$ and we might assume that the
fixed ordering of the roots is compatible with this decomposition. By Lemma \ref{determined} the factor $\v_1$ is completely determined by
the $c_\alpha$ for $\alpha\in \Delta'$ so there are finitely many 
possibilities for it. If there were infinitely many
elements in $\O\cap\w_0\v$, there would be infinitely many elements 
in $\O\cap\w_0\v_1 P^u$ for some $\v_1$.
We shall show that this cannot be the case. 
It is not restrictive to assume that
$G=SO_{2n+1}(k)$. We describe $G$ as the subgroup of $SL_{2n+1}(k)$ of matrices preserving the bilinear form associated with $\begin{pmatrix}1&0&0\cr
0&0&I_n\cr
0&I_n&0\cr\end{pmatrix}$ with respect to the canonical basis of $k^{2n+1}$. We may choose $B$ to be the subgroup of matrices of the form $\begin{pmatrix}1&0&^t\gamma\cr
-X\gamma&X&XA\cr
0&0&^tX^{-1}\cr\end{pmatrix}$ where $X$ is an invertible $n\times n$ upper triangular matrix, $^tX^{-1}$ is its inverse transpose, $\gamma$ is a column in $k^n$, $^t\gamma$ is its transpose and the symmetric part of $A$ is $-2^{-1}\gamma ^t\gamma$.
The above discussion and Lemma \ref{meet} translate into the assumption that there would be infinitely many conjugate matrices of the form 
\begin{align*}
x(V,\lambda)&=\w_0\v_1\v_2 =\begin{pmatrix}(-1)^n&0&0\cr
0&0&I\cr
0&I&0\cr\end{pmatrix}
\begin{pmatrix}1&0&0\cr
0&V&0\cr
0&0&^tV^{-1}\cr\end{pmatrix}
\begin{pmatrix}
1&0&^t\gamma\cr
-\gamma&I&A\cr
0&0&I\end{pmatrix}\\
&=\begin{pmatrix}(-1)^n&0&^t\!\gamma\cr
0&0&^t\!V^{-1}\cr
-V\gamma&V&VA\cr\end{pmatrix}
\end{align*}
where: $\gamma$ is a vector in $k^n$; $V$ is a fixed upper triangular unipotent matrix; $A$ is a matrix whose symmetric part is
$-2^{-1}\gamma ^t\!\gamma$ and by Lemma \ref{determined} the coefficients of $A$ and $\gamma$ depend polynomially on $\lambda=\gamma_n$ and the coefficients of $V$. 
The characteristic polynomial $q_\lambda(T)$ of $x(V,\lambda)$ depends polynomially on $\lambda$ thus $q_\lambda(T)=q_\mu(T)$ for at most finitely many $\mu$ in $k$
unless $q_\lambda(T)$ is independent of $\lambda$. We claim that this is not the case. In order to prove this, we need a more explicit description of $V$.\\
Using Lemma \ref{alpha-beta} one can show that, up to conjugation in $SO_{2n+1}(k)$ by diagonal matrices of type ${\rm diag}(1,-I_j,I_{n-j},-I_j,I_{n-j})$ the matrix $V$ is an upper triangular unipotent matrix with all $2$'s in the first off-diagonal. Inductively as in Lemma \ref{determined}, using Lemma \ref{alpha-beta} one sees that $V$ is the upper triangular unipotent matrix with only $2$'s above the diagonal. Thus it is enough to exhibit two matrices $x_1$ and $x_2$ of shape 
$\begin{pmatrix}(-1)^n&0&^t\!\gamma\cr
0&0&^t\!V^{-1}\cr
-V\gamma&V&M\cr\end{pmatrix}$ with $V$ as above,
lying in quasi-spherical conjugacy classes and with distinct characteristic polynomials.\\
Let $n$ be even. For $\zeta$ a square root of $2$ in $k$ we take $\gamma_1=2\zeta(1,-1,1,\,\cdots,\,-1)$ and 
$M_1=\begin{pmatrix}
0&-2&2&-2&\cdots\cr
2&-4&2&-2&\ddots\cr
2&-2&0&\ddots&\ddots\cr
\ddots&-2&\ddots&\ddots&\ddots\cr
2&\ddots&\ddots&\ddots&-4\cr
\end{pmatrix}
$. Then the matrix $x_1$ with $\gamma=\gamma_1$ and $M=M_1$ is conjugate to $a_1=\begin{pmatrix}
1&0&0\cr
0&-A_1&0\cr
0&0&-^t\!A_1\cr
\end{pmatrix}$ where $A_1$ is the upper triangular unipotent matrix with $(1,0,1,\,\cdots,0,1)$ on the first upper off-diagonal and zero elsewhere. Using the Jordan decomposition of $a_1$ and the formulas in \cite[\S 13.1]{carter-finite} which hold in good characteristic by \cite[Theorem 7.8]{Hu-cc} we see that the dimension of the conjugacy class $\O_{a_1}$ of $a_1$ is $n^2+n=\ell(w_0)+{\rm rk}(1-w_0)$. Since  $x_1$ lies in $\w_0 U$  we deduce from \cite[Theorem 5]{ccc}, \cite[Corollary 4.10]{gio} that $\O_{a_1}$ is spherical, hence quasi-spherical. Therefore $x=x_1(V,2\zeta)$. Let us now consider  $M_2=\begin{pmatrix}
4&-2&2&-2&\cdots\cr
2&0&2&-2&\ddots\cr
2&-2&4&\ddots&\ddots\cr
\ddots&-2&\ddots&\ddots&\ddots\cr
2&\ddots&\ddots&\ddots&0\cr
\end{pmatrix}
$. Then the matrix $x_2$ with $\gamma=0$ and $M=M_2$ is unipotent and lies in the conjugacy class corresponding to the Young diagram $(3,2^{n-2},1^2)$ whose dimension is again $n^2+n$. The class $\O_{x_2}$ is thus spherical (see also \cite[Theorem 3.2]{pany2} and \cite[Theorem 4.14]{FR}) hence quasi-spherical. It follows that $x_2=x(V,0)$ and the characteristic polynomials of $x_1$ and $x_2$ are different.\\ 
Let now $n$ be odd and let $\xi$ be a square root of $-2$ in $k$. We may consider $\gamma_3=-2\xi(1,-1,1,\cdots,-1,1)$ and 
$M_3$, constructed as $M_2$, 
and the corresponding matrix
$x_3$. One verifies that $x_3$ is unipotent and lies in the conjugacy class associated with the Young diagram $(3,2^{n-1})$, whose dimension is $n^2+n$. As above, this class is spherical, hence quasi-spherical so $x_3=x(V,-2\xi)$. On the other hand, taking $\gamma=0$ and $M_4$ constructed as $M_1$ 
we get a matrix $x_4$ which is conjugate to 
$a_2=\begin{pmatrix}
1&0&0\cr
0&-A_2&0\cr
0&0&-^t\!A_2\cr
\end{pmatrix}$ where $A_2$ is the upper triangular unipotent matrix with $(1,\,0,\,\cdots,1,\,\,0)$ on the first upper off-diagonal and zero elsewhere. As for $n$ even we see that the dimension of $\O_{a_2}$ is $n^2+n$ and since  $x_4$ lies in $\w_0 U$  we deduce as above that $\O_{a_2}$ is spherical, hence quasi-spherical. Thus, $x_4=x(V,0)$ and the characteristic polynomials of $x_3$ and $x_4$ are distinct. It follows that in a fixed class $\O$ there can only be finitely many elements of type $x(V,\lambda)$. By Lemma \ref{meet} each maximal $B$-orbit contains an element of type $x(V,\lambda)$ or a representative with all $c_\alpha=0$, so there are only finitely many of them.  
We may conclude using Lemma \ref{finite}.\hfill$\Box$

\begin{remark}\label{rsl2}{\rm  If $H$ is a connected reductive
    algebraic group the radical $R(H)$ of $H$ is a central torus (\cite[Proposition 7.3.1]{springer}) contained in all Borel subgroups of $H$. Thus, a conjugacy class
    $\O$ in $H$ is spherical (resp. quasi-spherical)  if and only if its projection into
   the semisimple group $K=H/R(H)$ is spherical (resp. quasi-spherical). Moreover, a conjugacy
    class in $K$ is spherical (resp. quasi-spherical) if and only if
    its projection into each simple factor of $K$ is spherical
    (resp. quasi-spherical). Thus, the results we obtained so far apply also to the case of reductive groups. In particular, if $G$ is connected, reductive with $w_0=-1$ and $\O$ is quasi-spherical in $G$ with $z_\O=w_0$ then $\O$ is spherical.}
\end{remark}

\section{Curves in the centralizer}\label{curves}

In this section we aim at the understanding of $G_x$ for $x\in \w
U\cap \O$ with $\O$ quasi-spherical and $w=z_\O$. We shall focus on
searching suitable families of elements in $G_x$ by making use
of the particular form of the chosen representative $x$ guaranteed by Lemma \ref{ortogonale}.  By Lemma \ref{U-} if
$w\alpha=\alpha$ then $X_{-\alpha}\subset X_\alpha s_\alpha B\cap G_x$. Now we aim at finding
elements in $X_{\gamma}s_\gamma B\cap G_x$ for the remaining roots $\gamma\in \Phi^+$, namely those such that $w\gamma\in
-\Phi^+$. We shall first
analyze those $\gamma$ for which $w\gamma=-\gamma$, that is,
$\gamma\in \Phi_1$, by looking at $G_x\cap G(\Phi_1)=G(\Phi_1)_x$. By Lemma \ref{ortogonale}, $x\in G(\Phi_1)$ so since
$w$ is the longest element in $W(\Phi_1)$ we may use
the results obtained in \S \ref{bigcell}. 

\begin{lemma}\label{dense-giuno}Let $\O$ be a quasi-spherical conjugacy class in
  $G$, let $w=z_\O$ and let $x\in \w U\cap\O$. Let $\alpha\in \Phi_1\cap\Delta$. 
Then for all but finitely  many $c\in k$ there exists $b_c\in B$ such that
$x_\alpha(c)n_\alpha b_c\in G_x\cap G(\Phi_1)$. 
\end{lemma}
\pf  We have $x\in G(\Phi_1)$ by
Lemma \ref{ortogonale}. The element $w$ is the longest element in $W(\Phi_1)$ and
its restriction to $W(\Phi_1)$ is $-1$. Let us consider the conjugacy class $\O'$  of $x$ in
$G(\Phi_1)$. It is quasi-spherical because $B_1=B\cap
G(\Phi_1)$ is a Borel subgroup of $G(\Phi_1)$ containing
$T$.
Therefore $\O'$ is spherical in $G(\Phi_1)$ by Remark \ref{rsl2}. 

Let $P$ be the minimal parabolic subgoup of $G(\Phi_1)$ associated with $\alpha$ and let $P^u$ be its unipotent radical. Let $x=\w x_{\alpha}(a)\v$ with $\v\in P^u$ and, for any nonzero $c\in k$, let
$y_c=n_\alpha^{-1}x_{\alpha}(-c)x x_{\alpha}(c)n_\alpha\in\O'$.
We have, for some nonzero structure constants
$\theta_1,\,\theta_2,\,\theta_3$ and for $t_1,\,t_2\in T$:
\begin{align*}
y_c&\in n_\alpha^{-1}(\w x_{-\alpha}(-\theta_1
c)x_\alpha(a+c)P^u) n_\alpha^{-1}\\
&= \w t_1  x_{\alpha}(-\theta_2\theta_1
c)x_{-\alpha}(\theta_3(a+c))P^u\\ 
&= \w t_1  x_{\alpha}(-\theta_2\theta_1
c)x_{-\alpha}((\theta_2\theta_1 c)^{-1})
x_{-\alpha}(-(\theta_2\theta_1 c)^{-1}+\theta_3(a+c))P^u\\
&= \w n_\alpha t_2 x_{\alpha}(\theta_2\theta_1
c) x_{-\alpha}(-(\theta_2\theta_1 c)^{-1}+\theta_3(a+c))P^u 
\end{align*}
where for the last equality we have used \cite[Lemma 8.1.4 (i)]{springer}.
Then, if $$-(\theta_2\theta_1 c)^{-1}+\theta_3(a+c)\neq 0$$ we have
$y_c\in B_1 w s_\alpha B_1s_\alpha B_1=B_1wB_1$ because $ws_\alpha < w$ in the Bruhat ordering. Thus,
for all but finitely many $c\in k$ the element $y_c$ lies in $B_1wB_1\cap
\O'$. By Lemma \ref{finite} applied to $\O'$ this intersection is the dense $B_1$-orbit $B_1.x$ so
for all but finitely many $c\in k$ there is $b_c\in B_1$ such that
$b_c^{-1} y_c b_c=x$, that is, $x_\alpha(c)n_\alpha b_c\in G_x\cap G(\Phi_1)$.
\hfill$\Box$

\bigskip

Next we shall consider the roots $\gamma$ for which
$w\gamma\neq\pm\gamma$. The set of all such roots is
$\Phi_2=\Phi\setminus(\Phi_1\cup\Phi(\Pi))$. Let $\Phi_2^+=\Phi_2\cap\Phi^+$. 
Every $\alpha\in\Phi_2^+$ determines
 the following subsets of $\Phi$:
\begin{equation*}\Phi^+(\alpha)=\bigcup_{j>0}(j\alpha+{\rm
    Ker}(1+w))\cap \Phi^+\subset\Phi^+(\Pi)\cup \Phi^+_2;
\end{equation*}
\begin{equation*}\Phi^-(\alpha)=\bigcup_{j>0}(j\alpha+{\rm
    Ker}(1+w))\cap (-\Phi^+)\subset (\Phi(\Pi)\cup \Phi_2)\cap(-\Phi^+)\end{equation*}
and $\Phi(\alpha)=\Phi^+(\alpha)\cup\Phi^-(\alpha)$. Fixing an
    ordering of the roots, we define:
\begin{equation}U_\alpha=\prod_{\beta\in\Phi^+(\alpha)}X_\beta\subset U;\quad\quad U_\alpha^-=\prod_{\beta\in\Phi^-(\alpha)}X_\beta\subset U^-.
\end{equation}

\begin{lemma}\label{lista}Let $\Phi$ be a simply or doubly-laced irreducible root system, let $w=w_0w_{\Pi}$ and let $\alpha\in\Phi^+_2$. We have:
\begin{enumerate}
\item The subsets $U_\alpha$ and $U^-_\alpha$ are subgroups of $U$
  and $U^-$, respectively.
\item The subgroups $U_\alpha$ and $U^-_\alpha$ are independent of the chosen ordering of the roots.
\item $U_\alpha\cap U_{\Phi_1}=1$.
\item If $X_{-\beta}\in U_{\alpha}^-$ then $X_{\beta}\not\in\langle U_{\Phi_1},\,U_\alpha\rangle=U_{\Phi_1}U_\alpha$.
\item $U_\alpha U^-_\alpha=U^-_\alpha U_\alpha$.
\item $w(\Phi(\alpha))\subset\Phi(\alpha)$
\item If $\Phi$ is simply-laced and
$w\alpha+\alpha\not\in-\Phi^+$ then $\w U_\alpha^- \w^{-1}\subset U_\alpha$.
\item If $\Phi$ is doubly-laced and $w\alpha+\alpha\not\in-(\Phi^+\cup 2\Phi^+)$ then $\w U_\alpha^- \w^{-1}\subset U_\alpha$.
\item If $\Phi$ is doubly-laced and $w\alpha+\alpha=2\beta\in -2\Phi^+$
  then $X_\beta\subset U_\alpha^-$ and
$\w U_\alpha^- \w^{-1}\subset U_\alpha X_\beta$. Besides, $X_\beta$ commutes with  $U_\alpha$.
\item If $\Phi$ is simply- or doubly-laced and $w\alpha+\alpha=\beta\in-\Phi^+$
  then $\w U_\alpha^- \w^{-1}\subset U_\alpha X_\beta$. Besides, $X_\beta$ commutes with $U_\alpha$.
\end{enumerate}
\end{lemma}
\pf The first two assertions follow from iterated application of \eqref{chev}. Statement 3 follows directly from the definition of $\Phi(\alpha)$.
Statement 4 is easily seen by looking at the coefficient of
  $\alpha$ in the expression of $\beta$. Statement 5 follows from 4 and \eqref{chev}. 
The sixth statement follows once we write $w=\prod_\gamma s_\gamma$
for mutually orthogonal roots $\gamma\in\Phi_1$.
Let us prove 7 and 8. If $\w U_\alpha^-\w^{-1}\subset U$ then $\w U_\alpha^-\w^{-1}\subset U_\alpha$ because $w$ is the product of reflections with respect to roots in $\Phi_1$. Hence, it is enough to show that $w\mu\in\Phi^+$ for all $\mu\in\Phi^-(\alpha)$. If we had $w\mu\in-\Phi^+$ for $\mu=j\alpha+y$ with $j>0$ and $y\in {\rm Ker}(1+w)$ we would have $\mu\in\Phi(\Pi)$ so $w\mu=\mu$, that is
\begin{equation}\label{deco}2\mu=\mu+w\mu=j\alpha+y+jw\alpha-y=j(\alpha+w\alpha)\in -2\Phi.\end{equation}
Thus $\alpha+w\alpha\not\in\Phi$ for it could neither be a positive nor a negative root, so $(\alpha,w\alpha)\geq0$. 
Taking $(2\mu,2\mu)$ we would have
\begin{equation}\label{gei}2(\mu,\mu)=j^2((\alpha,\alpha)+(\alpha,w\alpha))\geq
j^2(\alpha,\alpha).\end{equation}
If $(\alpha,\alpha)=(\mu,\mu)$ then $j=1$ and
$(\alpha,w\alpha)=(\alpha,\alpha)$ which is impossible proving statement 7. 
If $(\alpha,\alpha)=2(\mu,\mu)$ we have again $j=1$ and \eqref{deco} gives 
$2\mu=\alpha+w\alpha$ contradicting our assumption in the doubly-laced
case. 
If $2(\alpha,\alpha)=(\mu,\mu)$ 
we have $j^2\leq 4$ so $j\leq 2$.
Then either $j=2$ and
$\mu=\alpha+w\alpha\in-\Phi^+$ against our assumptions, 
or $j=1$ and
$3(\alpha,\alpha)=4(\alpha,w\alpha)$. Since this can never happen,
$\mu\not\in\Phi(\Pi)$  and statement 8 holds.\\
%
%
Let us prove 9. Let $\mu=j\alpha+y\in\Phi^-(\alpha)$, with $y\in{\rm Ker}(1+w)$
and $j>0$ and let us assume that $w\mu\in-\Phi^+$. It follows from the proof of 7 and 8 that we have $2\mu=j(\alpha+w\alpha)=2j\beta$. Hence $j=1$ and $\beta=\mu$ so $X_\beta$ is the only root subgroup in $U_\alpha^-$ that is mapped onto a negative root subgroup under conjugation by $\w$, and it is mapped
onto itself. Moreover, for every $\gamma=i\alpha+y'\in \Phi^+(\alpha)$ with $i>0$ and $y'\in{\rm Ker}(1+w)$ we have 
$2(\beta, \gamma)=(\alpha+w\alpha, i\alpha+y')=i(\alpha,\alpha)+i(\alpha,w\alpha)$ because $\alpha+w\alpha$ is orthogonal to ${\rm Ker}(1+w)$. Since $(\alpha,\alpha)=(w\alpha,w\alpha)$ we have $s_\alpha(w\alpha)\in\{w\alpha-\alpha,\, w\alpha,\,w\alpha+\alpha,\}$ so $2\frac{(\alpha,w\alpha)}{(\alpha,\alpha)}\in\{0,\pm1\}$. Thus $(\beta,\gamma)>0$ and therefore 
$\beta+\gamma\not\in\Phi$ so  $X_\beta$ commutes with with $X_\gamma$ and $\w U_\alpha^-\w^{-1}\subset U_\alpha X_\beta U_\alpha=U_\alpha X_\beta$. \\
Let us prove the last assertion. Let us assume that
 $\beta=\alpha+w\alpha\in-\Phi$. If for some root $\nu=j\alpha+y\in\Phi^-(\alpha)$
we had $w\nu\in-\Phi^+$ we would have, as before, $w\nu=\nu$ and
$2\nu=\nu+w\nu=j(\alpha+w\alpha)=j\beta\in 2\Phi$ so $j=2$ and $\beta=\nu$. Thus $\w U_\alpha^-\w^{-1}\subset U_\alpha X_\beta U_\alpha$. As in the proof of 9 we verify that $\beta+\gamma\not\in\Phi$ for every $\gamma\in\Phi^+(\alpha)$ whence $X_\beta$ commutes with $U_\alpha$ and $\w U_\alpha^-\w^{-1}\subset U_\alpha X_\beta$.\hfill$\Box$

\smallskip

\begin{lemma}\label{ra1}Let $G$ be a simple algebraic group, let $\O$ be a quasi-spherical conjugacy class with 
  $w=z_\O=w_0w_\Pi$ and let $x=\w\v\in\O\cap\w
  U$. Let
  $\alpha\in\Phi_2$ be such that $\alpha+w\alpha\not\in-\Phi^+$. Let us also assume, if $\Phi$ is
  doubly-laced, that $\alpha+w\alpha\not\in -2\Phi^+$.
Then   for every $c\in k$
  there exists an element in $x_{w\alpha}(c)U^w\cap G_x$.
\end{lemma}
\pf Since $w\alpha\neq\alpha$ we have $w\alpha\in-\Phi^+$.
For every $c\in k$ we consider the elements
\begin{align*}
y(c)&=x_{\alpha}(c)xx_{\alpha}(-c)=\w x_{w\alpha}(\theta c)\v x_{\alpha}(-c)\in\O
\end{align*}
where $\theta$ is a nonzero structure constant, and the elements
\begin{align*}
z(c)&=x_{w\alpha}(\theta c)\w x_{w\alpha}(\theta c)\v
x_{\alpha}(-c)x_{w\alpha}(-\theta c)\\
&=\w x_{\alpha}(\eta\theta c)(x_{w\alpha}(\theta c)\v x_{w\alpha}(-\theta c))
(x_{w\alpha}(\theta c)x_{\alpha}(-c) x_{w\alpha}(-\theta c))\in\O
\end{align*}
where $\eta$ is a nonzero structure constant. By making use of \eqref{chev} we shall show that for a
suitable $u_c\in U_\alpha$, possibly trivial, we
have
$$u^{-1}_cz(c)u_c=x\prod_{\gamma\in\Phi^+\setminus
  \Phi_1}x_\gamma(c_\gamma)\in\w U\cap \O.$$ Lemma \ref{ortogonale}
  will force
$c_\gamma=0$ for every $\gamma\in\Phi^+\setminus\Phi_1$ so $u^{-1}_cz(c)u_c=x$ and we have 
$u_c^{-1}x_{w\alpha}(\theta c)x_\alpha(c)=u'_cx_{\alpha}(c)x_{w\alpha}(\theta c)\in U_\alpha x_{w\alpha}(\theta c)\in G_x$. Taking inverses will give the statement because $c$ is arbitrary, $U_\alpha\subset U=U_wU^w$ and $U_w\subset G_x$ by Lemma \ref{ortogonale}.\\ 
By hypothesis $w\alpha+\alpha$ is either in $\Phi^+$ or it is not a root. Therefore we have $x_{w\alpha}(\theta c)x_{\alpha}(-c) x_{w\alpha}(-\theta c)=\v'\in U_\alpha$.
Besides, we have $\v=\prod_{\gamma_i\in\Phi^+_1}x_{\gamma_i}(c_i)$ so that
$$x_{w\alpha}(\theta c)\v x_{w\alpha}(-\theta
c)=\prod_{i=1}^r(x_{\gamma_i}(c_i)\prod_{a_i,b_i>0}x_{a_i\gamma_i+b_iw\alpha}(d_{abi}))$$
where we intend $d_{abi}=0$ if $a_i\gamma_i+b_iw\alpha\not\in\Phi$.
We proceed as follows: if  $w\alpha+\gamma_1\in-\Phi^+$ we apply
\eqref{chev} in order to move the term in $X_{w\alpha+\gamma_1}$ to the left of
$x_{\gamma_1}(c_1)$ whereas if $w\alpha+\gamma_1\in\Phi^+$ we apply
\eqref{chev} in order to move the term in $X_{w\alpha+\gamma_1}$ to the right of
$x_{\gamma_r}(c_r)$. At each step we might get extra factors either
in $U_\alpha$ or in $U^-_\alpha$ and we repeat the
procedure. Formula \eqref{chev} can always be applied
because we need never to interchange factors in $X_\beta$ with factors
in $X_{-\beta}$ (cfr. Lemma \ref{lista}(4)). Therefore we have:
$$x_{w\alpha}(\theta c)\v x_{w\alpha}(-\theta
c)=u^-\v u\in U_\alpha^-\v U_\alpha$$
because the coefficients of the terms
in $\v$ are never modified. By Lemma \ref{lista} (5) we have $x_\alpha(\eta\theta
c)u^-=u_-u_+\in U_\alpha^-U_\alpha$ and thus
$$z(c)=\w u_-u_+\v u\v'=\w u_-\v u'= u_c\w\v u'\subset U_\alpha \w\v U_\alpha$$
where for the second equality and the inclusions we have used Lemma \ref{lista} (4,7,8). 
Conjugation by $u_c^{-1}$ yields $u_c^{-1}z(c)u_c\in\O\cap \w\v U_\alpha$ hence
the term in $U_\alpha$ vanishes by Lemma \ref{ortogonale}. Then $u_c^{-1}x_{w\alpha}(\theta
c)x_\alpha(c)=u_c'x_\alpha(c)x_{w\alpha}(\theta c)\in G_x\cap
U_\alpha x_{w\alpha}(\theta c)$ and we have the statement. \hfill$\Box$
 
\smallskip

\begin{lemma}\label{w-pi-orbit}Let $\O$ be a quasi-spherical conjugacy
  class, let $w=z_\O=w_0w_\Pi$. If for some $\gamma\in \Phi^+\setminus\Phi(\Pi)$ and for every scalar $c$ there is 
  an element in $x_{-\gamma}(c)U^w$ centralizing  $x\in \O\cap \w U$ 
 then for every $\gamma'\in W_\Pi\gamma$ and for every $d\in k$ 
  there is an element in  $x_{-\gamma'}(d)U^w$ centralizing $x$.
\end{lemma}
\Pf By Lemmas \ref{ortogonale} and \ref{U-}, the centralizer of $x$
contains $X_{\pm \alpha}$ hence $n_\alpha$, for every
$\alpha\in\Pi$. Conjugation by $n_\alpha$ preserves $U^w$ and $U_w$ and maps
$X_{-\gamma}$ onto $X_{-s_\alpha(\gamma)}$, whence the statement.\hfill$\Box$

\smallskip

\begin{lemma}\label{duefi}Let $G$ be a simple algebraic group with  $\Phi$ doubly-laced. Let $\O$ be a
  quasi-spherical conjugacy class with notation as above. Let 
  $\alpha\in\Phi_2^+$ be such that 
  $w\alpha+\alpha=2\beta\in -2\Phi^+$. Then for $x\in\w U_\alpha\cap \O$ and for every
  $c\in k$ we have $x_{w\alpha}(c)U^w\cap G_x\neq\emptyset$.
\end{lemma}
\pf Let $z(c)$ be defined as in the proof of Lemma
\ref{ra1}. We have again  
$$x_{w\alpha}(\theta c)\v x_{w\alpha}(-\theta
c)=u^-\v u\in U_\alpha^-\v U_\alpha.$$
Let us first assume that 
$\beta\in-\Pi$. Then $u^-=x_{\beta}(a)u_-$ with $u_-\in U_\alpha^-\cap \w^{-1}U_\alpha
\w$ by Lemma \ref{lista} (2,9).
We have
$$
z(c)=\w x_\alpha(\theta\eta c) x_{\beta}(a)u_- \v u
x_{\alpha}(-c)=\w x_\beta(a)x_\alpha(\theta\eta c)u_-\v u x_\alpha(-c)$$
by Lemma \ref{lista} (9).  Applying repeatedly \eqref{chev} and Lemma \ref{lista} (5) we have for some $u'\in U_\alpha$, $u'_-,\,\v_-\in U^-_\alpha$, and $a'\in k$
$$z(c)=\w x_{\beta}(a)u_-'\v u'=\w x_{\beta}(a+a')\v_-\v u'$$
with $u_c=\w \v_-\w^{-1}\in U_\alpha$. 
We claim that $a+a'\neq0$. Otherwise, for some nonzero structure constant $\theta'$ we would have, 
by Lemma \ref{lista} (9),  
$$z(c)=x_{\beta}((a+a')\theta')u_c\w\v u'\in Bs_{\beta} Bw B=Bs_{\beta} wB$$
with $s_\beta w>w$ contradicting maximality of $w$. Thus, $a+a'=0$ and we
may proceed as in the proof of Lemma \ref{ra1}. Moreover, $u_c X_\alpha\subset U^w$.\\
If $\beta\not\in -\Pi$ then there is $\sigma\in W_{\Pi}$ such that $\sigma\beta\in -\Pi$ and $\sigma\alpha\in\Phi^+$  because the support of $\alpha$ contains at least one simple root outside $\Pi$. Since $w$ is the identity on $\Phi(\Pi)$ it commutes with $\sigma$ and we have $\sigma w\alpha\in-\Phi^+$ and $\sigma\alpha+w\sigma\alpha\in -2\Pi$. By the first part of the proof for every $c\in k$ there is an element in $x_{\sigma w\alpha}(c)U^w$ centralizing $x$ and we may apply Lemma \ref{w-pi-orbit} to get the statement.\hfill$\Box$

\smallskip

\begin{lemma}\label{vualfapiualfa}Let $G$ be a simple algebraic group with $\Phi$ simply or doubly-laced. Let $\O$
  be a quasi-spherical conjugacy class with notation as above. Let
  $\alpha\in\Phi^+$ be such that  $w\alpha+\alpha=\beta\in-\Phi$. Then for $x\in\w U\cap \O$ and for every
  $c\in k$ we have $x_{w\alpha}(c)U^w\cap G_x\neq\emptyset$.
\end{lemma}
\pf Let us first assume that $\beta\in -\Pi$.
%
We have $X_\beta\subset
U_\alpha^-$ and $\w U^-_\alpha \w^{-1}\subset X_\beta U_\alpha$ by Lemma \ref{lista}(10).
We have $\beta+w\alpha=2w\alpha+\alpha=w(w\alpha+2\alpha)\not\in\Phi$. This follows as in the proof of Lemma \ref{lista} (9,10).

Let $z(c)$ be as in the as in the proofs of Lemmas \ref{ra1} and \ref{duefi}. As above we have:
$$x_{w\alpha}(\theta c)\v x_{w\alpha}(-\theta
c)=u^-\v u\in U_\alpha^-\v U_\alpha.$$
We may apply \eqref{chev} to move $x_{\alpha}(\eta\theta c)$ to the right of
$\v$ in the expression of $z(c)$. Then we have
$$z(c)=\w u_-\v u^+ x_{w\alpha}(\theta c)x_\alpha(-c)x_{w\alpha}(-\theta c)$$
with $u_-=x_\beta(h)u_-'\in X_\beta U_\alpha^-$, $u_c=\w u_-'\w^{-1}\in U_\alpha\cap U^w$ and $u^+\in U_\alpha$. Applying once more \eqref{chev} to
$ x_{w\alpha}(\theta c)x_\alpha(-c)$ gives only a nontrivial extra term
in $X_{\beta}$ by Lemma \ref{lista} (10).
%
Then, for some $h_1,\,h_2\in k$ and some $u'\in U_\alpha$ we have
$z(c)=x_{\beta}(h_1)u_c\w\v u' x_\beta(h_2)$. Conjugation by
$u_c^{-1}x_{\beta}(-h_1)$ yields an element $z'(c)$ in $x U_\alpha
X_{\beta}U_\alpha\cap\O\subset BwBX_\beta B\subset Bws_\beta B\cup BwB$. Maximality of $w$ forces $h_1+h_2=0$ so
the $X_\beta$-factor in $z'(c)$ is trivial. Lemma
\ref{ortogonale} implies that $z'(c)=x$ so, using that $\beta+\alpha,\,\beta+w\alpha\not\in\Phi$ we have
$$u_c^{-1}x_\beta(-h_1)x_{w\alpha}(\theta
c)x_\alpha(c)=u_c^{-1}x_\alpha(c)x_{w\alpha}(\theta c)x_\beta(h)\in
G_x$$ with $u_c\in U^w$. By Lemma \ref{U-} we have $x_\beta(h)\in G_x$ so  $u_c^{-1}x_\alpha(c)x_{w\alpha}(\theta c)\in
G_x$ and taking the inverse yields the statement for $\beta\in-\Pi$.\\
If $\beta\not\in -\Pi$ we may apply Lemma \ref{w-pi-orbit} as we did in Lemma \ref{duefi}.
\hfill$\Box$
\smallskip

We have constructed enough elements in $G_x$ and we are ready to prove the main result of this paper.

\begin{theorem}\label{quasiquasi} Let $G$ be a connected, reductive algebraic group over an algebraically closed field $k$ of characteristic zero or good and odd. Then
every quasi-spherical conjugacy class $\O$ in $G$ is spherical.
\end{theorem}
\pf By Remark \ref{rsl2} it is enough to prove the statement for $G$ simple. 
Type $G_2$ has already been discussed in Section \ref{g2} so we
only need to consider $\Phi$ simply or doubly-laced. Moreover, when $z_\O=w_0=-1$ the statement has been proved in Propositions \ref{bigcellDE}, \ref{bigcellf4}, \ref{cn}, \ref{bn} so we shall prove the remaining cases. 
Let $v$ be a maximal $B$-orbit in $\O$.
We will prove the statement by showing that $\dim(\O)=\dim(v)$,
so that $v$ is dense in $\O$. To this end, we need to show that for
some $x\in v$ we have
$\dim G_x=\dim B_x+|\Phi^+|$. We will do so by using $x\in\w U\cap\O$ for $w=z_\O$.

Let us consider the restriction $\pi_x$ to $G_x$ of the natural
projection $\pi$ of $G$ onto the flag variety $G/B$. Let $gB$ be in
the image of $\pi_x$. We may assume that $g\in G_x$ and then it is not hard
to verify that $\pi_x^{-1}(gB)=gB_x$ so each non-empty fiber has dimension equal to $\dim B_x$. Since
$\dim G/B=|\Phi^+|$ it is
enough to prove that $\pi_x$ is dominant  and use \cite[Theorem 5.1.6]{springer}. We shall prove that
$\pi_x(G_x)\cap \pi(B\sigma B)$ is dense in $\pi(B\sigma B)$ for every
$\sigma\in W$. In particular, this is true for $\sigma=w_0$ so
$\pi_x(G_x)\cap \pi(Bw_0B)$ is dense in $\pi(Bw_0 B)$ thus
$\pi_x(G_x)$ is dense in $G/B$. 

More precisely, if we identify $\pi(B\sigma B)=\pi(U^\sigma
\sigma B)$ with the affine space ${\mathbb A}^{\ell(\sigma)}$ through
the map $\pi(u\sigmad B)=\pi(\prod_{\gamma\in\Phi_\sigma}x_\gamma(c_\gamma)\sigmad B)\mapsto
(c_\gamma)_{\gamma\in\Phi_\sigma}$, we will show by induction on $\ell(\sigma)$ that
$\pi_x(G_x)\cap {\mathbb A}^{\ell(\sigma)}$ contains the complement in
${\mathbb A}^{\ell(\sigma)}$ of finitely many hyperplanes. 

For $\sigma=1$ there is nothing to say. Suppose that
the statement holds for $\ell(\sigma)\leq s$ and let us consider
$\tau\in W$ with $\ell(\tau)=s+1$. Then $\tau=\sigma s_\alpha$ for
some $\sigma\in W$ with $\ell(\sigma)=s$ and some $\alpha\in\Delta$ with $\sigma\alpha\in\Phi^+$. 
 Besides, $\Phi_\tau=\Phi_\sigma\cup\{\sigma\alpha\}$ so $U^\tau=U^\sigma X_{\sigma\alpha}$. 
 By the induction hypothesis the set $U'$ of elements $u$ in $U^\sigma$ for which
$u\sigmad b$ lies in $G_x$ for some $b\in B$ contains the complement of finitely many hyperplanes in $U^\sigma\cong {\mathbb A}^{\ell(\sigma)}$.
 
There are three possibilities: $\alpha\in\Pi$, $\alpha\in\Delta\cap\Phi_1$
and $\alpha\in\Delta\cap\Phi_2$.\\ 
If $\alpha\in\Pi$ we have $X_{\alpha}n_\alpha\subset G_x$ by Lemma \ref{ortogonale}. 
Then for every $u\in U'$ and every $c\in k$ there is $b\in B$ for which $(u\sigmad b)(x_{\alpha}(c)n_\alpha)\in G_x$. Let $b=x_{\alpha}(r)\v$ for $r\in k$ and $\v\in P^u_\alpha$. Then for some $\v'\in P_\alpha^u$ and for some nonzero structure constant $\eta$ we have
$$(u\sigmad b)(x_{\alpha}(c)n_\alpha)=u \sigmad x_{\alpha}(r+c)n_\alpha \v'=u x_{\sigma\alpha}(\eta(r+c))\sigmad n_\alpha \v'\in G_x.$$
Since $c$ is arbitrary and $\eta\neq0$, if $\alpha\in\Pi$ then 
$\pi_x(G_x)\cap \pi(B\tau B)$ contains $\pi(U'X_{\sigma\alpha}\tau B)$ so it contains the complement of finitely many hyperplanes in ${\mathbb A}^{\ell(\tau)}$. 

Let now $\alpha\in\Delta\cap\Phi_1$. By Lemma \ref{dense-giuno} for all but finitely many $c\in k$ there is $b_c\in B$ such that $x_\alpha(c)n_\alpha b_c\in G_x$. Thus, for every $u\in U'$ and for those $c$ there is $b\in B$ for which $(u\sigmad b)(x_{\alpha}(c)n_\alpha b_c)\in G_x$. Let $b=x_{\alpha}(r)\v$ for $r\in k$ and $\v\in P^u_\alpha$. Then for some $\v'\in P_\alpha^u$ and for some nonzero structure constant $\eta$ we have
$$(u\sigmad b)(x_{\alpha}(c)n_\alpha b_c)=u \sigmad x_{\alpha}(r+c)n_\alpha \v' b_c=u x_{\sigma\alpha}(\eta(r+c))\sigmad n_\alpha \v'b_c\in G_x.$$
Since all but finitely many $c$ were possible and $\eta\neq0$, also in this case $\pi_x(G_x)\cap \pi(B\tau B)$ contains $\pi(U'x_{\sigma\alpha}(c)\tau B)$ for all but finitely many $c$, thus it contains the complement of finitely many hyperplanes in ${\mathbb A}^{\ell(\tau)}$. 

Finally, let $\alpha\in\Phi_2\cap\Delta$. Then by Lemmas \ref{ra1},
\ref{duefi}, and \ref{vualfapiualfa} for every $c\in k$ there exists
$u_c\in U$ such that $x_{-\alpha}(c)u_c\in G_x$. For $c\neq0$ this element is
equal to $x_{\alpha}(c^{-1})t_cn_\alpha
x_{\alpha}(c^{-1})u_c$ for some $t_c\in T$ by \cite[Lemma 8.1.4 (i)]{springer}. Thus, for
every $u\in U'$ and $c\neq0$ there is $b\in B$ for which $(u\sigmad b)(x_{\alpha}(c^{-1})t_c n_\alpha
x_{\alpha}(c^{-1})u_c)\in G_x$. Let $b=x_{\alpha}(r)\v$ for $r\in k$
and $\v\in P^u_\alpha$. Then for some $\v'\in P_\alpha^u$, $t_c'\in T$ and for some nonzero structure constant $\eta$ we have
\begin{align*}
(u\sigmad b)(x_{\alpha}(c^{-1})t_c n_\alpha
x_{\alpha}(c^{-1})u_c)
&=u \sigmad x_{\alpha}(r+c^{-1})n_\alpha t_c'\v' x_{\alpha}(c^{-1})u_c\\
&=u x_{\sigma\alpha}(\eta(r+c^{-1}))\sigmad n_\alpha t_c'\v' x_\alpha(c^{-1}) u_c\in G_x.
\end{align*}
Then again, $\pi_x(G_x)\cap \pi(B\tau B)$ contains the complement of finitely many hyperplanes in ${\mathbb A}^{\ell(\tau)}$ and we have the statement.\hfill$\Box$

\smallskip

As a consequence of Theorem \ref{quasiquasi} we get the sought characterization.

\begin{theorem}\label{cara}Let $\O$ be a conjugacy class in a connected reductive algebraic
  group $G$ over a field of zero or good odd characteristic. Then
  $\O$ is spherical if and only if $\O\subset \bigcup_{w^2=1} BwB$.
\end{theorem}
\pf This is obtained combining Theorem \ref{quasiquasi} with
\cite[Theorem 2.7]{gio}, whose proof holds also for $G$ connected and reductive.\hfill$\Box$

\section*{Acknowledgements}

The author is grateful to the referee for suggesting improvements in the presentation of the manuscript and to Mauro Costantini for useful comments on a previous version of this paper.

\end{document}